\documentclass[preprint, 10pt]{elsarticle}

\usepackage[margin=2.5cm]{geometry}
\usepackage{setspace}
\usepackage{lmodern}
\usepackage{amsmath,amssymb}
\usepackage{amsfonts,amsthm}
\usepackage{lineno}
\usepackage{graphicx}
\usepackage{bm, bbm}
\usepackage{algorithm}
\usepackage{physics}
\usepackage{color}
\usepackage{pxfonts}
\usepackage[footnotesize,bf]{caption}
\usepackage{subcaption}
\usepackage{mathtools}
\usepackage{hhline}
\usepackage[T1]{fontenc}

\usepackage{caption}
\usepackage{graphicx}
\usepackage{float} 
\usepackage{subcaption}

\usepackage{placeins}

\newcommand{\R}{\mathbb{R}}

\newtheorem{remark}{Remark}[section]

\makeatletter
\def\ps@pprintTitle{
	\let\@oddhead\@empty
	\let\@evenhead\@empty
	\let\@oddfoot\@empty
	\let\@evenfoot\@oddfoot
}
\makeatother

\graphicspath{{figure/}}
\usepackage{multirow}
\usepackage{ulem}
\usepackage{tikz}
\usepackage{color}
\usepackage[]{xcolor}

\newcounter{mycommand}
\renewcommand{\themycommand}{S\arabic{mycommand}}
\newcommand{\mycommand}{\refstepcounter{mycommand}\themycommand}

\usepackage{titlesec}
\titleformat{\section}{\bf}{S\thesection}{}{}

\usepackage{makecell}

\biboptions{sort&compress}

\begin{document}
\begin{frontmatter}
\title{IB-UQ: Information bottleneck based uncertainty quantification for neural function regression and neural operator learning}
\author[SNH]{Ling Guo}
\author[SJTU,TJ]{Hao Wu\corref{cor}}
\author[SNH]{Wenwen Zhou}
\author[TJ]{Yan Wang}
\author[LESC]{Tao Zhou}
\date{}
\address[SNH]{Department of Mathematics, Shanghai Normal University, Shanghai, China}
\address[SJTU]{School of Mathematical Sciences, Institute of Natural Sciences, and MOE-LSC, Shanghai Jiaotong University, Shanghai, China}
\address[TJ]{School of Mathematical Sciences, Tongji University, Shanghai, China}
\address[LESC]{LSEC, Institute of Computational Mathematics and Scientific/Engineering
Computing, Academy of Mathematics and Systems Science, Chinese Academy
of Sciences, Beijing, China.}
\cortext[cor]{Corresponding Author}

\begin{abstract}
\iffalse
We propose a novel framework  for uncertainty quantification via information bottleneck (IB-UQ) for scientific machine learning tasks, including deep neural network (DNN) regression and neural operator learning (DeepONet). Specifically, we first employ the General Incompressible-Flow Networks (GIN) model to learn a "wide" distribution from noisy observation data. Then, following the information bottleneck objective, we learn a stochastic map from input to some latent representation that can be used to predict the output. A tractable variational bound on the IB objective is constructed with
a normalizing flow reparameterization.  Hence, we can optimize the objective using the stochastic gradient descent method. IB-UQ  can provide both mean and variance in the label prediction by explicitly modeling the representation variables. Compared to most DNN regression methods and the deterministic DeepONet, the proposed model can be trained on noisy data and provide accurate predictions with reliable uncertainty estimates on unseen noisy data. 
We demonstrate the capability of the proposed IB-UQ framework via several representative examples, including discontinuous function regression, real-world dataset regression and learning nonlinear operators for diffusion-reaction partial differential equation.
\fi

We propose a novel framework for uncertainty quantification via information bottleneck (IB-UQ) for scientific machine learning tasks, including deep neural network (DNN) regression and neural operator learning (DeepONet). Specifically, we
incorporate the bottleneck by a confidence-aware encoder, which encodes inputs into latent representations according to the confidence of the input data belonging to the region where training data is located,
and utilize a Gaussian decoder to predict means and variances of outputs conditional on representation variables. Furthermore, we propose a data augmentation based  information bottleneck objective which can enhance the quantification quality of the extrapolation uncertainty, and the encoder and decoder can be both trained by minimizing a tractable variational bound of the objective.
In comparison to uncertainty quantification (UQ) methods for scientific learning tasks that rely on Bayesian neural networks with Hamiltonian Monte Carlo posterior estimators, the model we propose is computationally efficient, particularly when dealing with large-scale data sets. The effectiveness of the IB-UQ model has been demonstrated through several representative examples, such as regression for discontinuous functions, real-world data set regression, learning nonlinear operators for partial differential equations, and a large-scale climate model. The experimental results indicate that the IB-UQ model can handle noisy data, generate robust predictions, and provide confident uncertainty evaluation for out-of-distribution data.

\end{abstract}

\end{frontmatter}

\section*{\bf\LARGE{Main text}}

Scientific machine learning, particularly physics-informed deep learning, has achieved significant success
in modelling and predicting the response of complex physical systems \cite{karniadakis2021physics,willard2021integrating}. Deep neural networks (DNNs), which possess the universal approximation property for continuous functions \cite{chen1993approximations},  has been widely employed in approximating the solutions of forward and inverse ordinary/partial differential equations \cite{Lagaris1997,khoo2021solving,MaziarParisGK17_1,LiaoMing2021,guo2022normalizing,guo2022monte,huang2022augmented,gao2022failure}. They have also been proven useful in solving high-dimensional partial differential
equations \cite{EYu2018,ZangBaoYeZhou2020} and discovering equations from data \cite{brunton2016discovering,long2018pde}.   Another prominent direction in scientific machine learning lie in learning continuous operators or complex systems from streams of scattered data. One representitave approach is the deep operator network (DeepONet) proposed by Lu et al. \cite{lu2021learning}. DeepONet's original structure is based on the universal approximation theorem of nonlinear operators \cite{chen1995universal} and has been successfully used in multiscale and multiphysics problems \cite{lin2021operator,mao2021deepm}. To further enhance the performance of DeepONet, a physics-informed architecture was proposed in \cite{wang2021learning}. Another structure for operator learning is the Fourier Neural Operator \cite{li2020fourier}, which parameterizes the integral kernel in Fourier space. Recently, a novel kernel-coupled attention operator learning method was developed in \cite{kissas2022learning} inspired by the attention mechanism's success in deep learning.

There are several sources of uncertainty that can arise when using DNNs in scientific machine learning tasks, including noisy and limited data, neural network parameters, unknown parameters, and incomplete physical models, among others \cite{psaros2022uncertainty}. Thus it is crucial to accurately quantify and predict uncertainty for deep learning to be reliably used in practical applications. 
 Proper uncertainty estimates are essential to evaluate a model's confidence in its predictions. The most popular approach  of uncertainty quantification (UQ) in the deep learning community is based on the Bayesian framework \cite{mackay1995bayesian,neal2012bayesian,gal2016dropout}, where the uncertainty is obtained by performing posterior inference using Bayes’ rule given observational data and prior beliefs. Alternative UQ methods are based on ensembles of DNN optimization iterates or independently trained DNNs \cite{lakshminarayanan2017simple,malinin2018predictive,fort2019deep}, as well as on the evidential framework \cite{sensoy2018evidential,amini2020deep}. Uncertainty quantification (UQ) in the context of scientific machine learning is a more challenging task due to the inclusion of physical models. There have been few works on UQ for scientific machine learning especially for operator learning till now. Psaros et al. provided an extensive review and proposed novel methods of UQ for scientific machine learning in \cite{psaros2022uncertainty}, accompanied by an open-source Python library termed NeuralUQ \cite{zou2022neuraluq}. Other related UQ research for operator learning and reliable DeepONet method can be found in \cite{moya2023deeponet,lin2021accelerated,yang2022scalable,zhu2022reliable}. 

The information bottleneck (IB) principle was initially introduced by Tishby et al. in their work \cite{tishby1999information}. The principle seeks to find an encoding $Z$ that can maximally express the target variable $Y$ while being maximally
compressive about the input $X$. To achieve this, the principle proposes maximizing the objective function $I(Z;Y)-\beta I(Z;X)$, where $\beta$ is the Lagrange multiplier. The theoretical framework of the IB principle is used to analyze deep neural networks in \cite{tishby2015deep}. A variational approximation of IB (VIB) was then established in \cite{alemi2016deep}, where a nerual network is used to parameterise the information bottleneck model. Further results in \cite{alemi2018uncertainty} show VIB incorporates uncertainty naturally thus can improve calibration and detect out-of-distribution data. Addational related research including different approximation approaches for the IB objective, the connection between IB and the stochastic gradient descent training dynamics can be found in \cite{kolchinsky2019nonlinear,shwartz2017opening,saxe2019information} and references therein. While many existing studies on IB-based uncertainty quantification have focused on classification tasks, there has been limited attention given to regression tasks thus far.

\begin{figure}[htbp]
	\centering
	\includegraphics[width=0.98\linewidth]{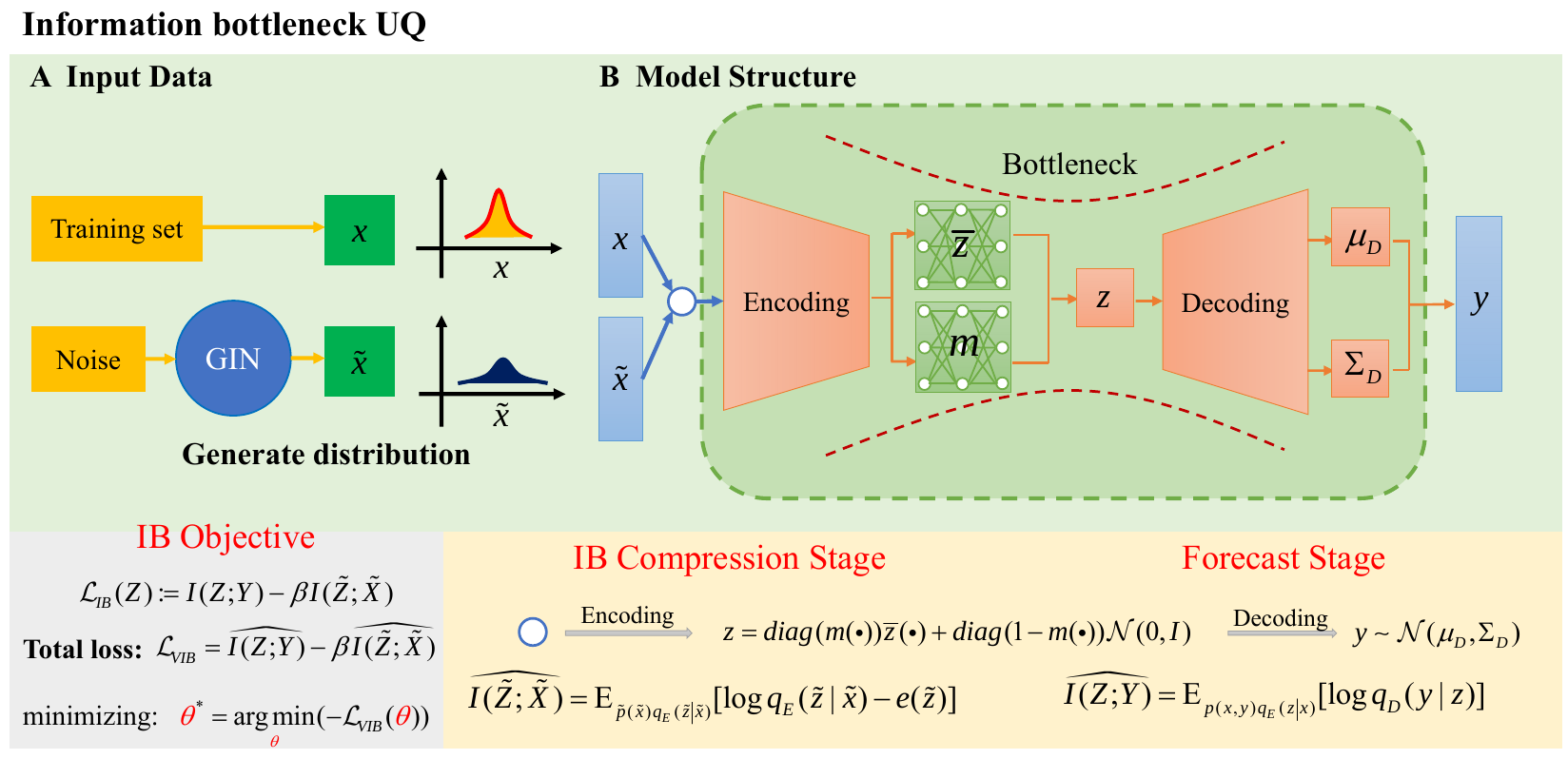}
	\caption{IB-UQ: Schematic of information bottleneck based uncertainty quantification for function regression with deep neural networks. 
    The input is $x$ taken from the training data or $\tilde x$ distributed according to a wide distribution given by a General Incompressible-Flow (GIN). After passing an encoder, the input is encoded as a random variable $z$ in the latent space by extracting features $\bar z(x)$ and incorporating noise with indensity $1-m(x)$. Finally, a Gaussian conditional distribution of the output can be obtained by a decoder, where both the mean $\mu_D$ and covariance $\Sigma_D$ are functions of $z$. Here, $\bar z, m, \mu_D, \Sigma_D$ are all modeled by neural networks, and can be trained by maximizing a variational lower bound $\mathcal L_{\mathrm{VIB}}$ of the IB objective $\mathcal L_\mathrm{IB}$. }
	\label{fig:IBUQ_fun}
\end{figure}

In this paper, we put effort into establishing a novel framework for uncertainty quantification via information bottleneck (IB-UQ) in scientific machine learning tasks, including deep neural network regression and DeepONet operator learning. A schematic of IB-UQ method for function regression is shown in Figure \ref{fig:IBUQ_fun}, which is specifically designed to improve uncertainty prediction for both in-distribution and out-of-distribution inputs. To achieve this, we employ the general incompressible-flow network (GIN) to generate augmented data from a "wide" distribution, that has a high entropy and can cover the domain of the training data, and
establish a data augmentation based IB objective. Additionally, we introduce a confidence-aware encoder that maps the input to a latent representation that can be used to predict the mean and covariance of the output, allowing for more accurate and reliable uncertainty estimates. To ensure that our model can be trained efficiently, we construct a tractable variational bound on the IB objective using a normalizing flow reparameterization, which enables all parameters of the model to be trained according to the bound. The main contributions of this work can be summarized as follows: 

\begin{itemize}
\item  We introduce a novel framework for quantifying uncertainty in scientific machine learning tasks, including deep neural network (DNN) regression and neural operator learning (DeepONet).
\item  The proposed IB-UQ method can obtain both mean and uncertainty estimates of the prediction via explicitly modeling the random representation variables. In comparison to the existing gold standard UQ method, such as Bayesian Neural Network with Hamiltonian Monte Carlo posterior estimation \cite{psaros2022uncertainty}, the proposed new model is both easy to implement and efficient for large-scale data problems.

\item The proposed model is capable of providing confident uncertainty evaluation for out-of-distribution data, which is crucial for the use of neural network-based methods in risk-related tasks and applications. This capability enhances the model's reliability and applicability in real-world scenarios where the presence of such data is inevitable.
\end{itemize}

\section*{\bf\LARGE{Results}}

We demonstrate the effectiveness of our proposed IB-UQ method on various deep neural network (DNN) function regression and deep operator learning problems. 
We will present detailed experimental results for several benchmark problems, including regression of a discontinuous function, operator learning for two partial differential equations, California housing prices regression, and a large-scale climate modeling task.  Additionally, we report more experimental analysis and comparison results for the discontinuous function regression problem and climate modeling in \ref{app:additional-results}.   The hyperparameters and neural network architectures we used in these examples are summarized in Supplementary Information \ref{app:h-param}.

The numerical results indicate the IB-UQ model's ability to provide robust functional predictions, while also demonstrate its ability to handle noisy data and provide confident uncertainty estimates for out-of-distribution  samples.

The codes used to generate the results will be
released in GitHub upon publication of the paper.

\subsection*{\bf\large{Discontinuous function regression}}
We first consider a one-dimensional discontinuous test function 
\begin{equation}\label{eq:DG function}
u(x)=\left\{\begin{array}{l}
\frac{1}{2}\left[\sin ^3(2 \pi x)-1\right],-1 \leq x<0, \\
\frac{1}{2}\left[\sin ^3(3 \pi x)+1\right], 0 \leq x \leq 1.
\end{array}\right.
\end{equation}
The training dataset $\mathcal{D}=\{ x_{i},u_{i} \}^{N}_{i=1}$ consists of $N=32$, unless otherwise specified, equidistant measurements of $u(x)$ at $x_i \in [-0.8, -0.2] \cup [0.2, 0.8]$. The data is contaminated with zero-mean Gaussian noise $\epsilon_u\sim \mathcal{N}(0, \sigma_u^{2})$.

\begin{figure}[htp!]
\centering
\includegraphics[width=0.75\linewidth]{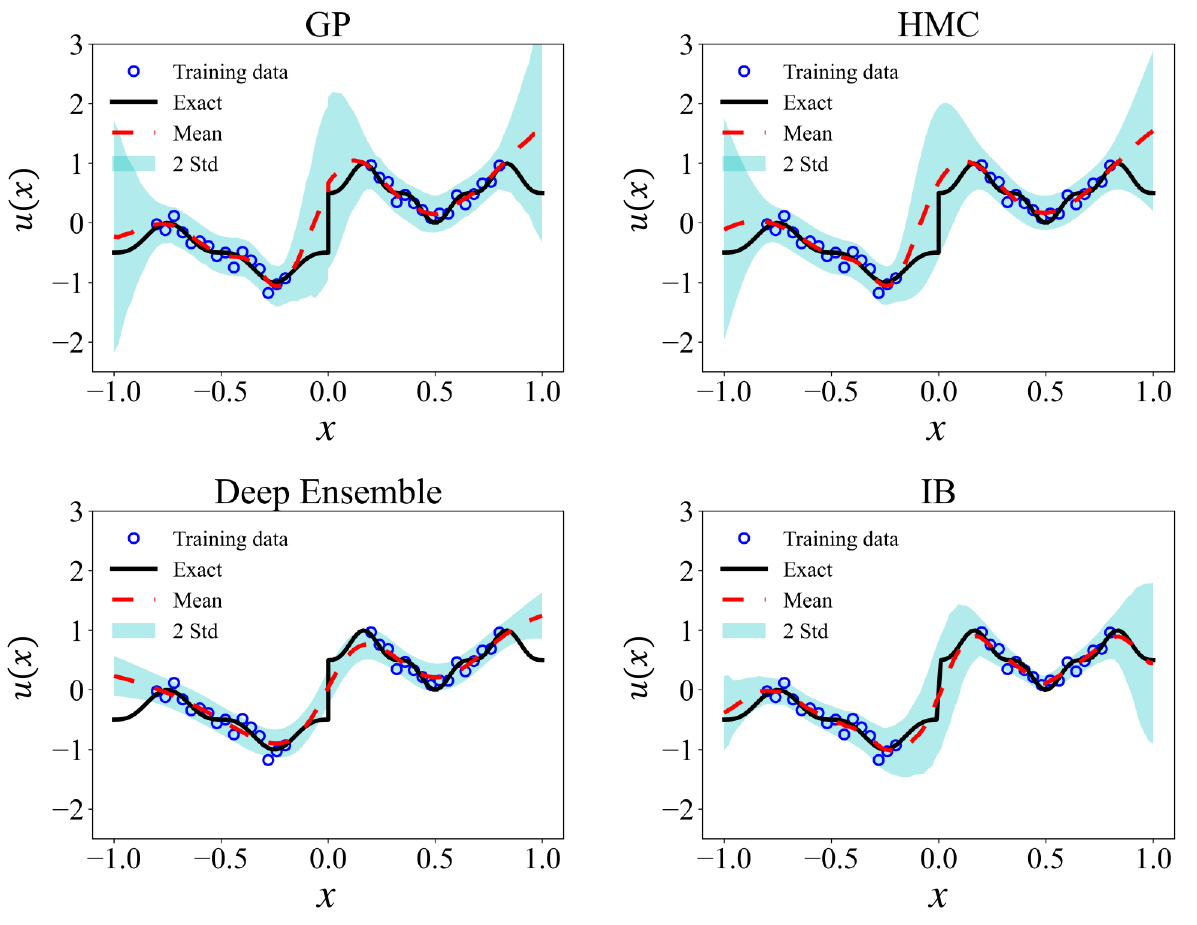}
\caption{Function regression problem (\ref{eq:DG function}) with noise scale $\sigma_u=0.1$: comparison among different approaches. The training data and exact function, as well as the mean and uncertainty (the standard deviations) are shown here. We can see that all the methods can get mean prediction that is similar to the exact solution while  the uncertainty estimates cover in most cases the point-wise errors. }
\label{fig:fun_pre}
\end{figure}

We use the IB-UQ algorithm to predict the values of $u(x)$ at any unseen position $x$ with uncertainty
estimates with $\beta=0.3$, and compare IB-UQ method with three efficient UQ methods for scientific machine learning problems given in \cite{psaros2022uncertainty}, including Gaussian processes regression (GP), Bayesian Neural networks with Hamiltonian Monte Carlo posterior sampling method  (HMC), and Deep Ensemble method. Figure \ref{fig:fun_pre} presents the mean prediction and uncertainty estimates, i.e. the standard deviations, obtained by the four UQ methods. We can see that (1) all the four methods provide good mean prediction compared with the exact solution; (2) similar as GP and HMC methods, IB-UQ method can also provide larger standard deviations (i.e., uncertainty) at the regions with fewer training data, i.e., $x\in [-1, -0.8]\bigcup [-0.2, 0.2]\bigcup [0.8, 1]$, which is important in risk-sensitive applications, where computational methods should be able to distinguish between interpolation domain and extrapolation domain. However, the uncertainty of the Deep Ensemble method does not increase obviously for the out-of-distribution (OOD) test data. Although GP and HMC are often considered the gold standard for Bayesian neural network (BNN) methods, they both have limitations. Gaussian processes regression, for instance, is hard to cope with nonlinearities when applied to solve partial differential equations (PDEs). Additionally, the HMC method can be computationally expensive, particularly for large data sets, compared to our IB-UQ method, which does not require posterior distribution sampling and is easy to implement.

\subsection*{\bf\large{Diffusion-reaction equation}}

Here we consider the following diffusion-reaction equation 

\begin{equation}\label{diffusion reaction equation}
 \quad\quad\quad  \frac{\partial s}{\partial t}=D \frac{\partial^2 s}{\partial x^2}+k s^2+u(x), \quad x \in [0,1], \quad t \in [0,1],
\end{equation}
with zero initial and boundary conditions \cite{zhu2022reliable}. Here $D$ denotes the diffusion coefficient and $k$ the reaction rate. In this example, $D$ is set at $0.01$ and $k$ is set at $0.5$. We aim to use our information bottleneck based UQ method to learn the operator that maps the source term $u(x)$ to the solution $s(x; t)$ of the system (\ref{diffusion reaction equation}) with : 
$$
\mathcal{G}: u(x) \longmapsto s(x,t). 
$$

The training and testing data are generated as follows. The input functions $u(x)$ are sampled from a mean-zero Gaussian random process $u(x)\sim \mathcal{GP}(0, k_l(x_1,x_2))$ with an squared-exponential kernel  $$k_l(x_1,x_2)=\text{exp}\bigg (- \frac{\|x_1-x_2\|^2}{2l^2}\bigg),$$
where $l > 0$ represents the correlation length. The training dataset uses $l=0.5$ and each realization of $u(x)$ is recorded on an equi-spaced grid of $m = 101$ sensor locations. Then we solve the diffusion-reaction system using a second-order implicit finite-difference method on a $101 \times 101$ grid to obtain the training data of $s(x,t)$.

We totally generate $N = 10000$ input/output function pairs for training the IB-UQ model. To report the performance of the IB-UQ model for out-of-distribution data detection, the $1000$ test data are obtained from Gaussian process with different correlation length scales. Specifically, for in-distribution (ID) test data, $u(x)$ corresponds to a sample from the same stochastic process as used for training, i.e $l=0.5$. For out-of-distribution (OOD) test data, $u(x)$ corresponds to a sample from a stochastic process with different correlation length $l$ of the training data, i.e. $l> 0.5$ or $l< 0.5$. We consider the scenario that both the training/test data values are contaminated with Gaussian noise with standard deviation equal to $0.01$. The regularization parameter $\beta$ in IB-UQ objective function \eqref{eq:obj-aug} for operator learning is set to $\beta=0.3$.

\begin{figure}[!htp]
\centering
\includegraphics[width=0.7\linewidth]{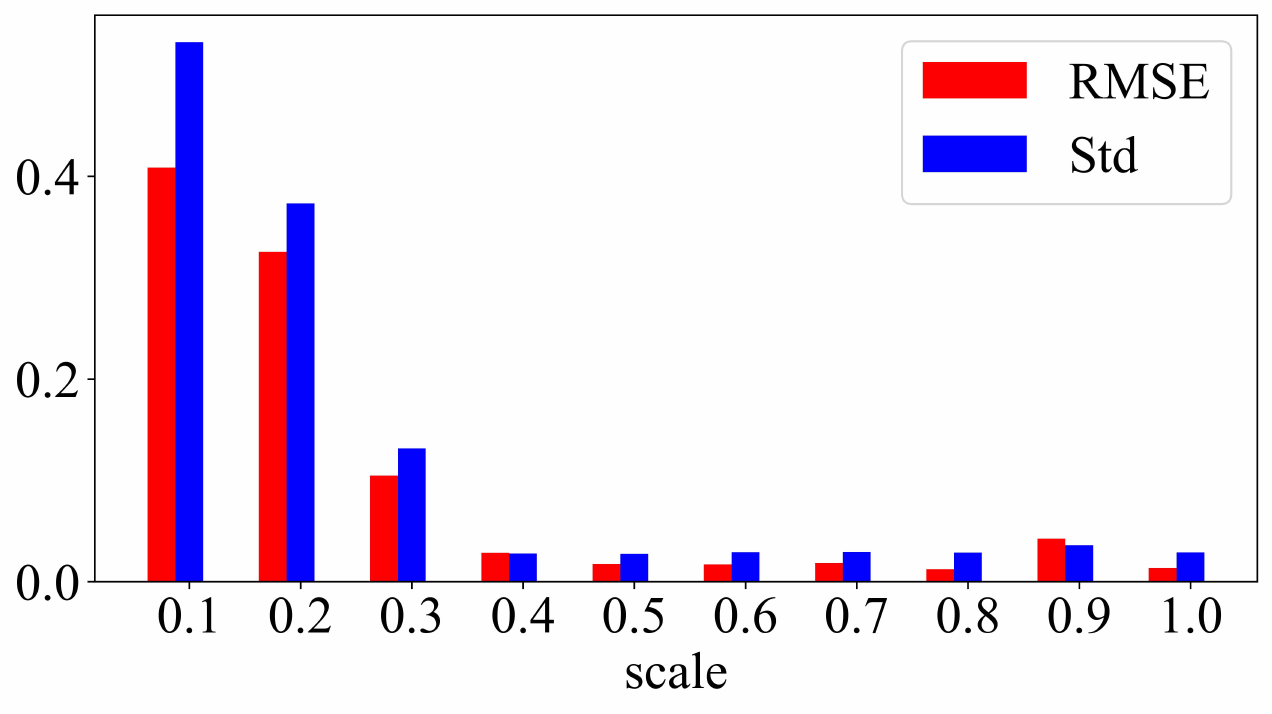}
\caption{Operator learning of Reaction-diffusion equation. Shown here are the RMSE values of the IB-UQ mean prediction for inputs $u$ associated with the different correlation lengths $l$ and the associated standard deviation (1 std) predicted by IB-UQ. All values are averaged over 1000 inputs.}
\label{fig:Onetmrse}
\end{figure}

\begin{figure}[htp!]
\centering
\subfloat{\includegraphics[width=0.85\linewidth]{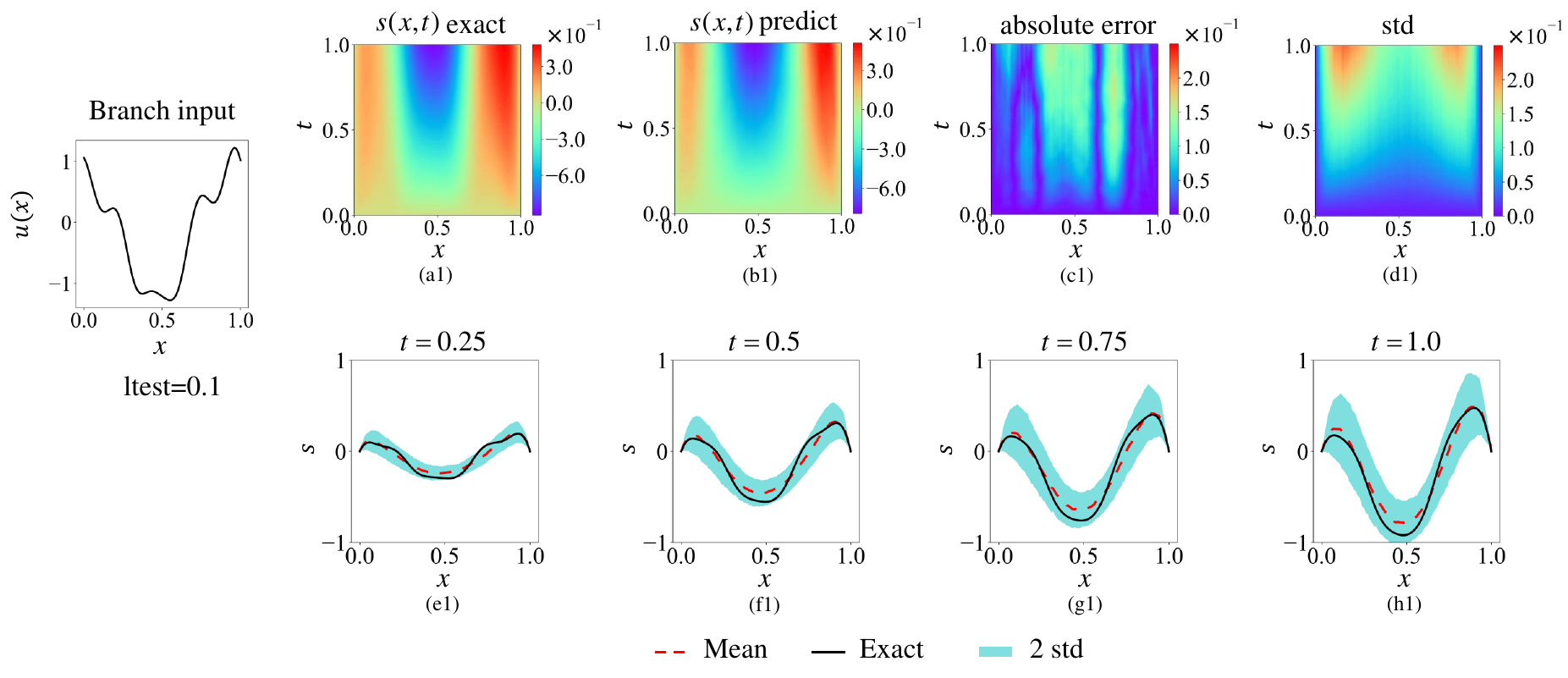}}\\
\subfloat{\includegraphics[width=0.85\linewidth]{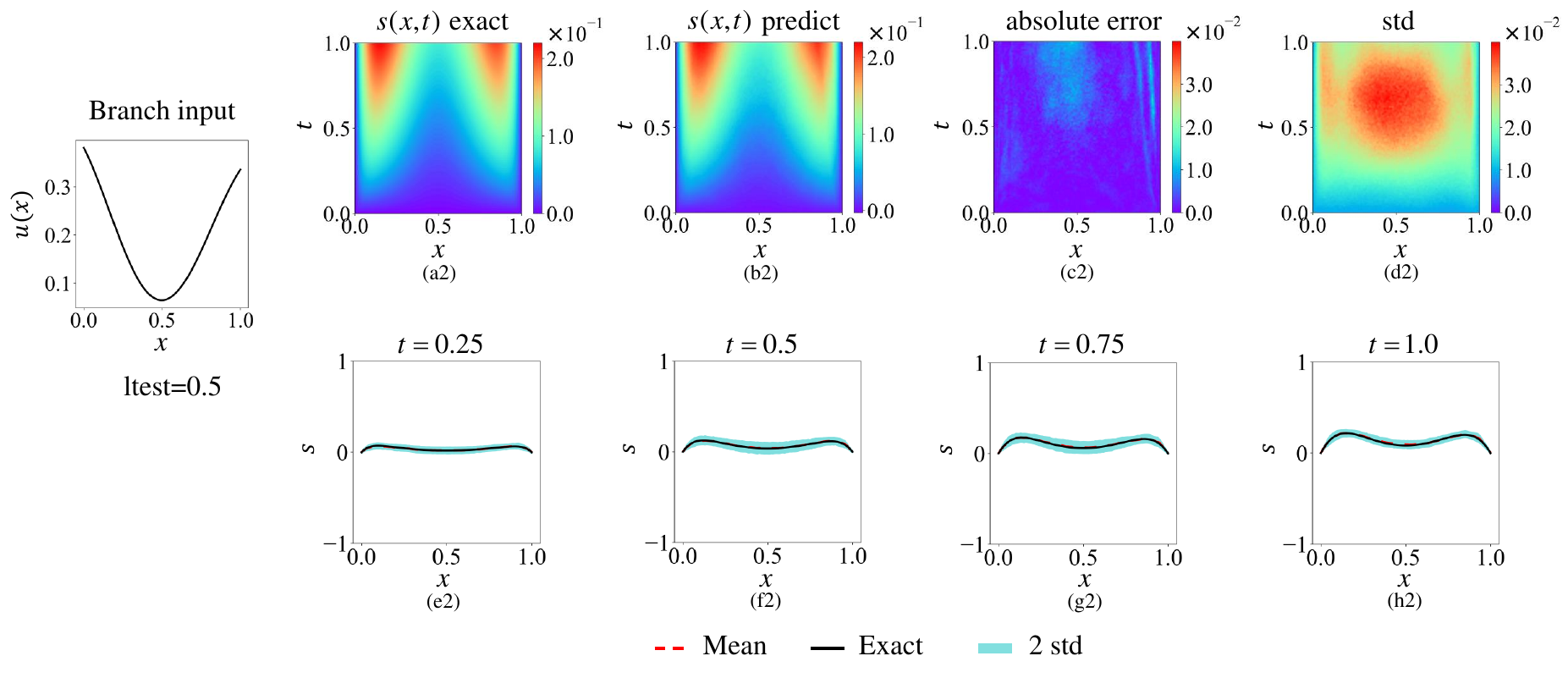}}\\
\subfloat{\includegraphics[width=0.85\linewidth]{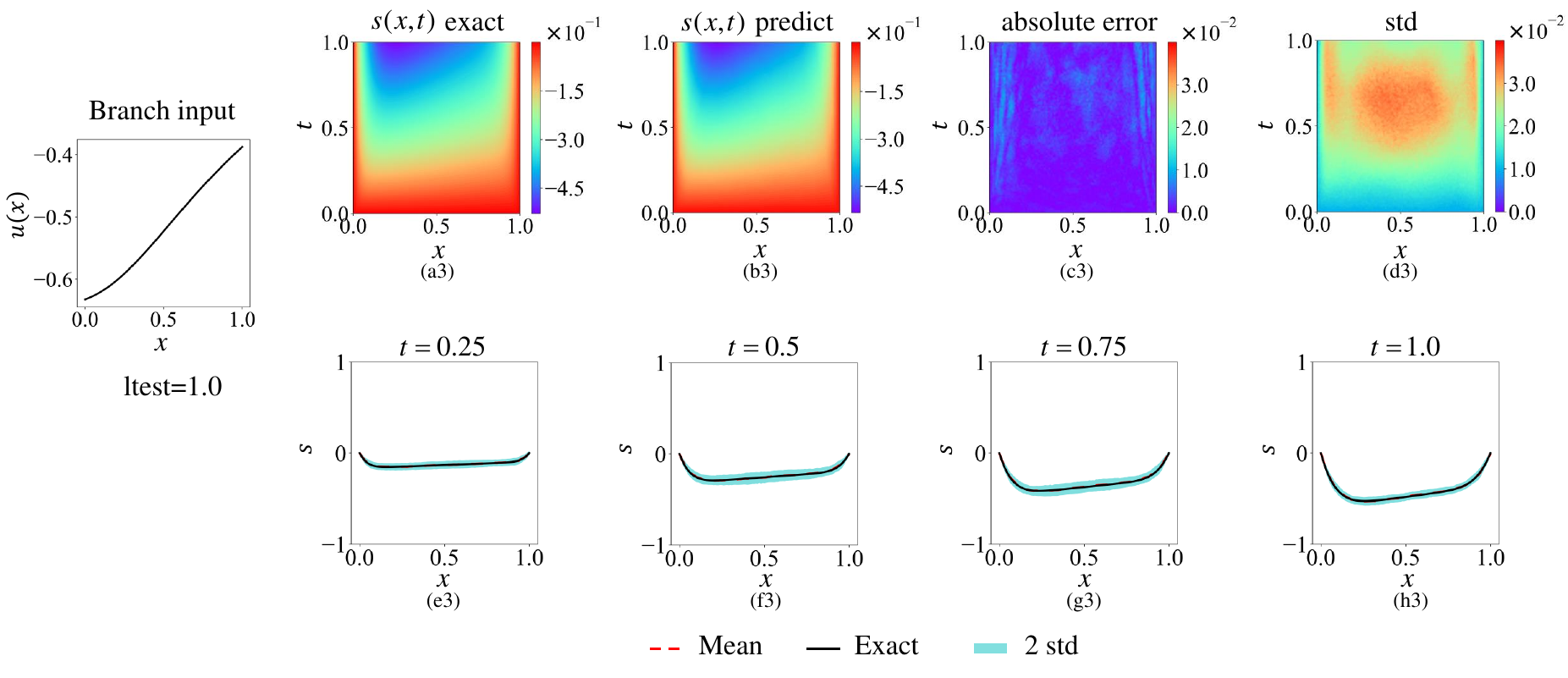}}
\caption{Operator learning of reaction-diffusion equation. Representative results of IB-UQ method with different input test functions (left frame). {\bf{OOD test results with}} $\mathbf{l=0.1}$: 
 (a1)-(d1) Ground truth solution,  IB-UQ DeepONet predictive mean, the absolute error between  the ground truth solution and the IB-UQ predictive mean, as well as 
 the standard deviations (predictive uncertainty); (e1)-(h1) Ground truth solution (black), IB-UQ predictive mean (red dashed) and uncertainty (blue shade) at t = 0.25, t = 0.50, t = 0.75, and t=1 respectively. {\bf{ID test results with}} $\mathbf{l=0.5}$:  (a2)-(h2) show the IB-UQ output visualization corresponding to test data with input functions having the same correlation length as the training data.  {\bf{OOD test results with}} $\mathbf{l=1.0}$: (a3)-(h3) plot the IB-UQ output visualization corresponding to OOD test data, where the input function is sampled from a Gaussian process with correlation length larger than the training set.}
\label{fig:dr}
\end{figure}

The root mean square error (RMSE) between the ground truth solution and the IB-UQ DeepONet predictive mean as well as the uncertainty estimate are reported in Figure  \ref{fig:Onetmrse}. We can see that when the training and test functions are sampled from the same space ($l=0.5$), i.e. ID test data,  the error is smaller. However, the RMSE of IB-UQ model increases significantly for OOD data with $l< 0.5$. Nevertheless, for OOD data with $l> 0.5$, IB-UQ DeepONet still has a small error. Potential explanation is with a larger $l$ leading to smoother kernel functions and DeepONet structure can predict accurately for smoother functions. These results are consistent with the observations reported in \cite{lu2021learning,yang2022scalable}. Figure \ref{fig:dr} presents comparison between the ground truth solution and the IB-UQ predictive mean and uncertainty tested on data generated with different correlation length. We obsever that the computational absolute errors between the predicted means and the reference solutions are bounded by the predicted uncertainties, 2 standard deviation, for the the three different cases. Moreover, despite the predictions of the IB-UQ not being accurate on testing data with smaller $l=0.1$, the larger uncertainty obtained from IB-UQ reflects the inaccuracy of the predictions. Thus IB-UQ model can provide reliable uncertainty estimates for OOD data. 

%%%%%%%%%%%%%%%%%%%%%%%%%%%%%%%%%%%%%%%%%%%%%%%%%%%%%%%%%
%%%%%%%%%%%%%%%%%%%%%%%%%%%%%%%%%%%%%%%%%%%%%%%%%%%%%%%%%

\subsection*{\bf\large{Advection equation}}

We now consider to train the IB-UQ model to learn the operator that  maps $v(x)$ to the solution $u(x,t)$:
$$
\mathcal{G}: v(x) \longmapsto u(x,t),
$$
defined by the following advection equation \cite{zhu2022reliable}
\begin{equation}\label{advection equation}
    \frac{\partial u}{\partial t} + v(x) \frac{\partial u}{\partial x} = 0, \quad x \in [0,1], \quad t \in [0,1],
\end{equation}
with the initial conditional $u(x,0)=\text{sin}(\pi x)$ and boundary conditioal $u(0,t)=\text{sin}(\pi t /2)$. We take $v(x)$ in form of $v(x)=V(x)-\text{min}_{x}V(x)+1$ to guarantee $v(x)>0$. 

\begin{figure}[htp!]
\centering
\subfloat{\includegraphics[width=0.85\linewidth]{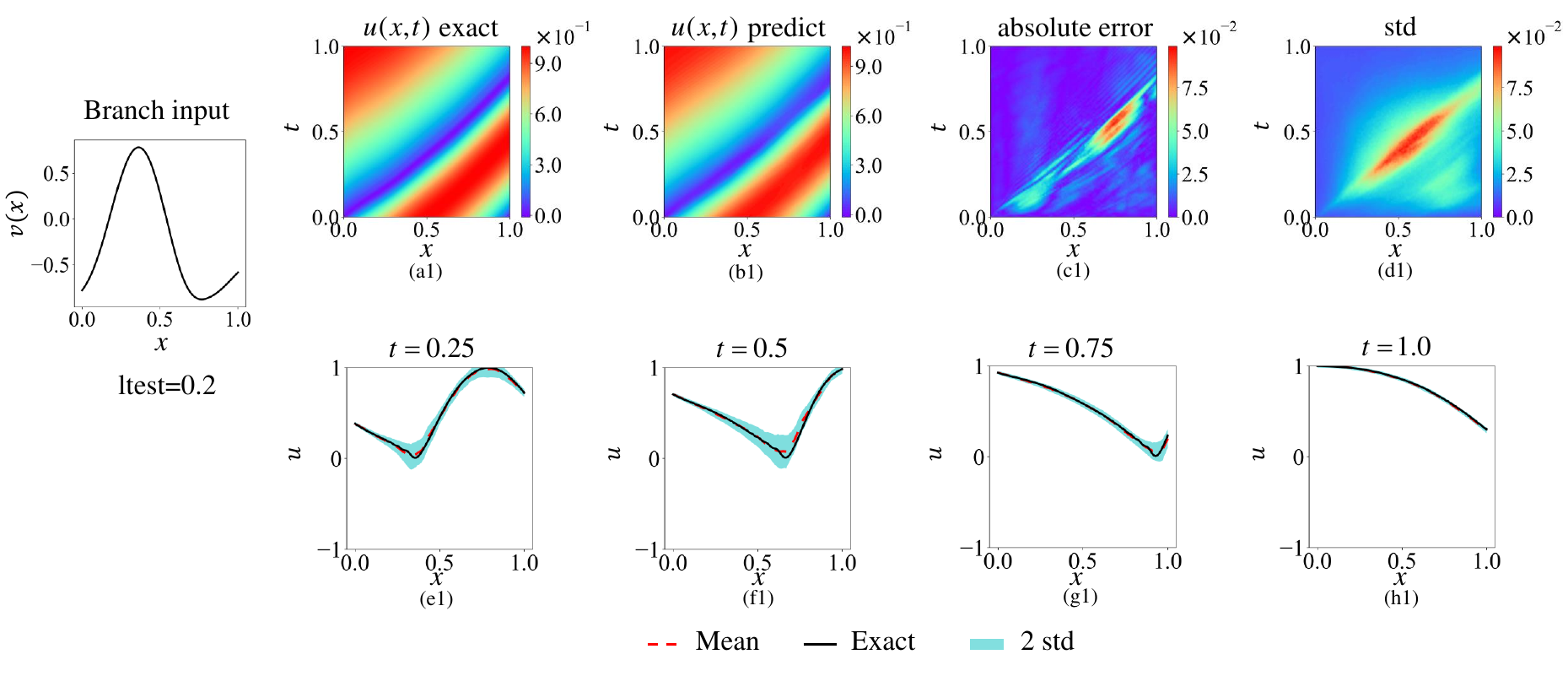}}\\
\subfloat{\includegraphics[width=0.85\linewidth]{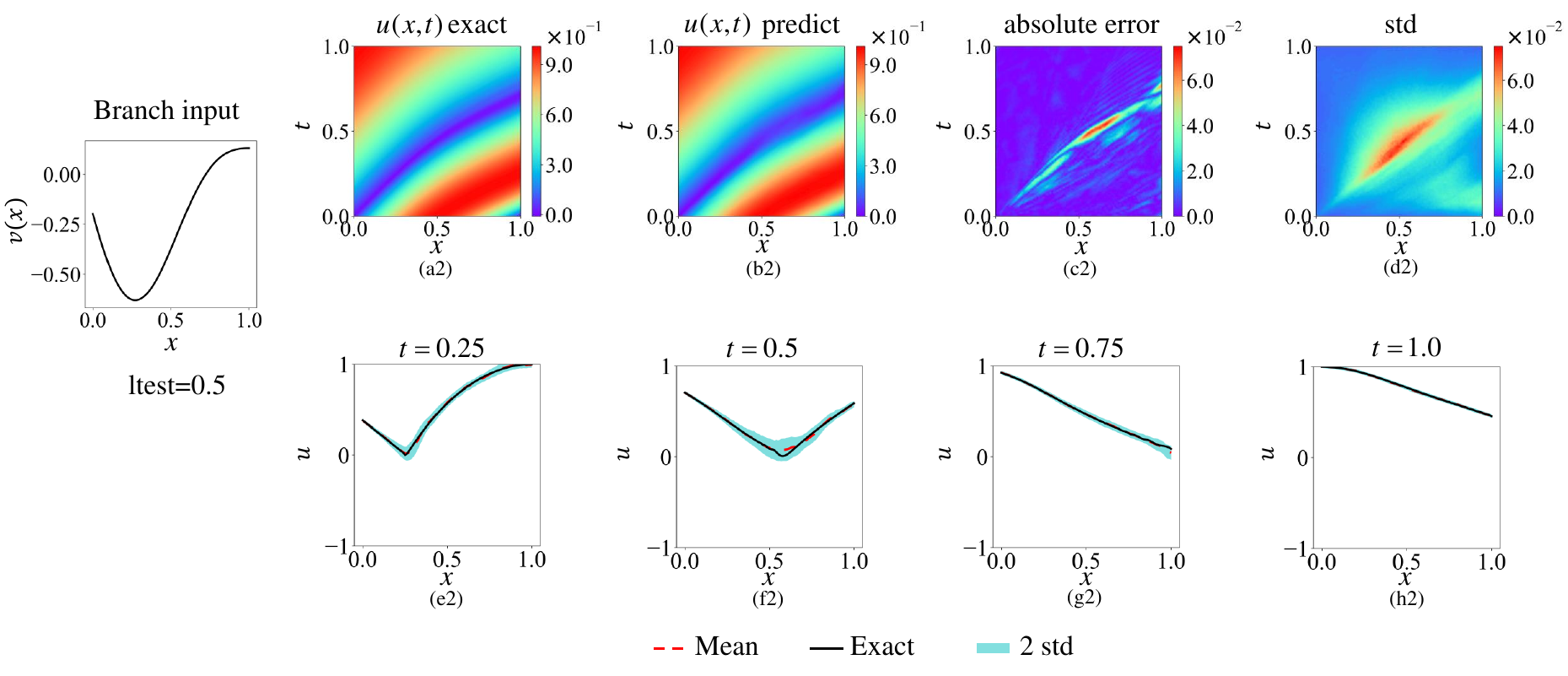}}\\
\subfloat{\includegraphics[width=0.85\linewidth]{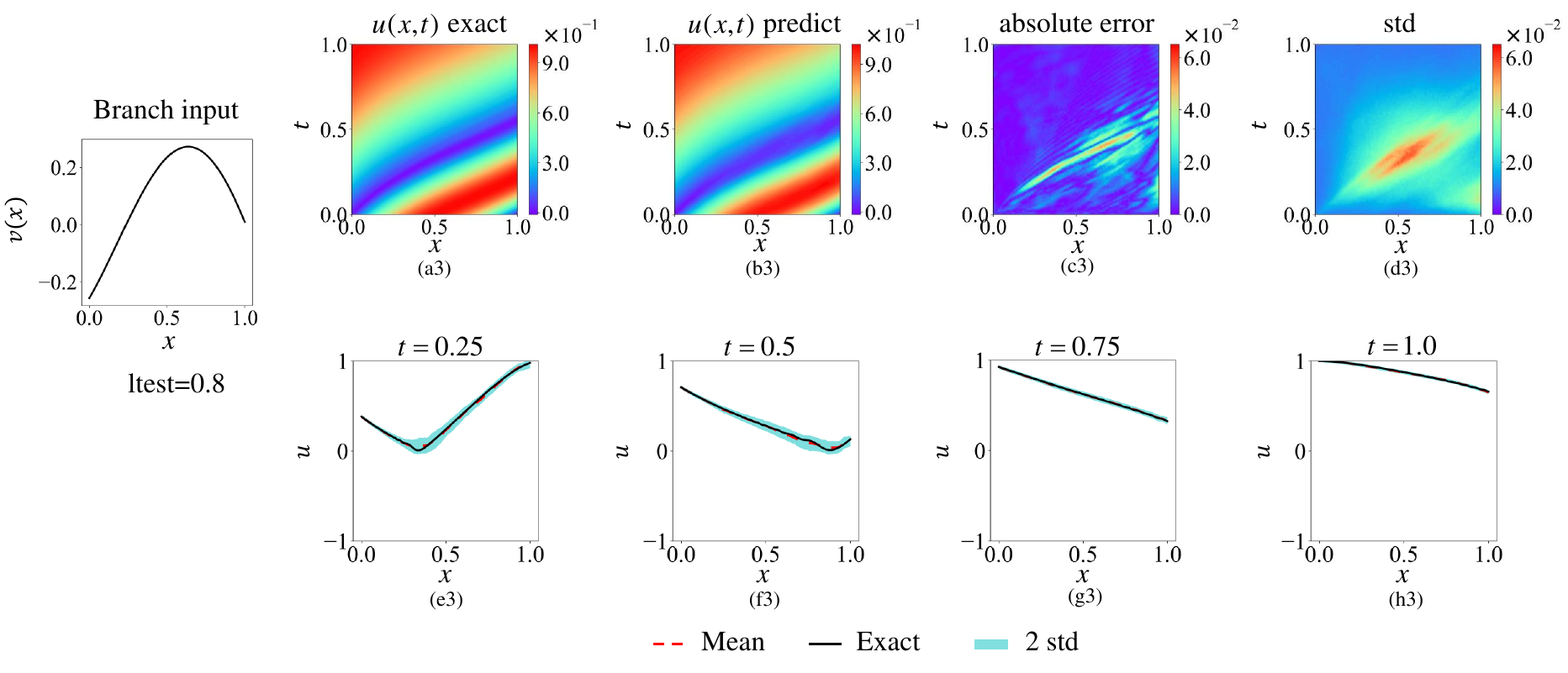}}
\caption{Operator learning of advection equation. Representative results of IB-UQ method with different input test functions (left frame). {\bf{OOD test results with}} $\mathbf{l=0.2}$: 
 (a1)-(d1) Ground truth solution,  IB-UQ DeepONet predictive mean, the absolute error between  the ground truth solution and the IB-UQ predictive mean, as well as 
 the standard deviations (predictive uncertainty); (e1)-(h1) Ground truth solution (black), IB-UQ predictive mean (red dashed) and uncertainty (blue shade) at t = 0.25, t = 0.50, t = 0.75, and t=1 respectively. {\bf{ID test results with}} $\mathbf{l=0.5}$:  (a2)-(h2) show the IB-UQ output visualization corresponding to test data with input functions having the same correlation length as the training data.  {\bf{OOD test results with}} $\mathbf{l=0.8}$: (a3)-(h3) plot the IB-UQ output visualization corresponding to OOD test data, where the input function is sampled from a Gaussian process with correlation length larger than the training set.}
	\label{fig:adv} 
\end{figure}

The training data set is constructed by sampling $v(x)$ from a Gaussian random field with an squared exponential kernel with correlation length $l=0.5$, and then use them as inputs to solve the advection equation via a finite difference method to obtain $N=10000$ input and output function pairs, where the measurements on $101\times 101$ grid points of each $u(x,t)$ are available. To test the performance of the IB-UQ method on ID and OOD data, we generate three test datasets of 1000 functions with $l_{test} = 0.2,0.5,0.8$ respectively.
We consider the scenario that both training/test data values are contaminated with Gaussian noise with standard deviation equal to 0.01. 

We present the comparison between the ground truth
solution and the IB-UQ predictive mean and uncertainty in Figure \ref{fig:adv} for testing cases with different correlation length $l$. We also plot 4 time-snapshots from 1 out of 1000 randomly choosing testing samples. We can see that the absolute errors between the predicted means and the reference solutions are bounded by the predicted uncertainties. We also observe that for OOD test sample generated with correlation length $l=0.2$, the prediction mean is less accurate than those of ID and OOD test sample with larger correlation length $l=0.2$. Thus larger uncertainty is obtained in the cases where the predictive mean is not accurate. This example also demonstrate that the IB-UQ model can provide resonable predictions with confident uncertainty estimates on noisy unseen dataset. 

%%%%%%%%%%%%%%%%%%%%%%%%%%%%%%%%%%%%%%%%%%%%%%%%%%%%
%%%%%%%%%%%%%%%%%%%%%%%%%%%%%%%%%%%%%%%%%%%%%%%%%%%%%

\subsection*{{\bf\large{California housing prices data set}}}

\begin{figure}[!htp]
\centering
\subfloat[]{\includegraphics[width=0.33\linewidth]{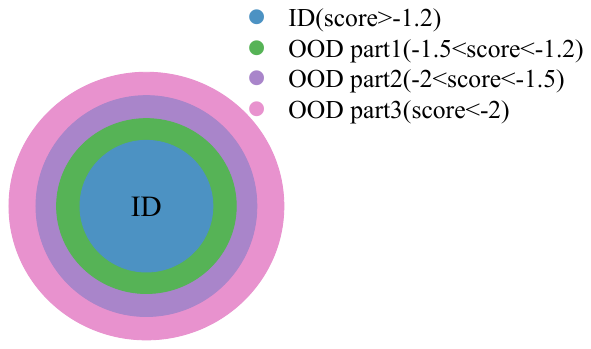}}\hfill
\subfloat[]{\includegraphics[width=0.33\linewidth]{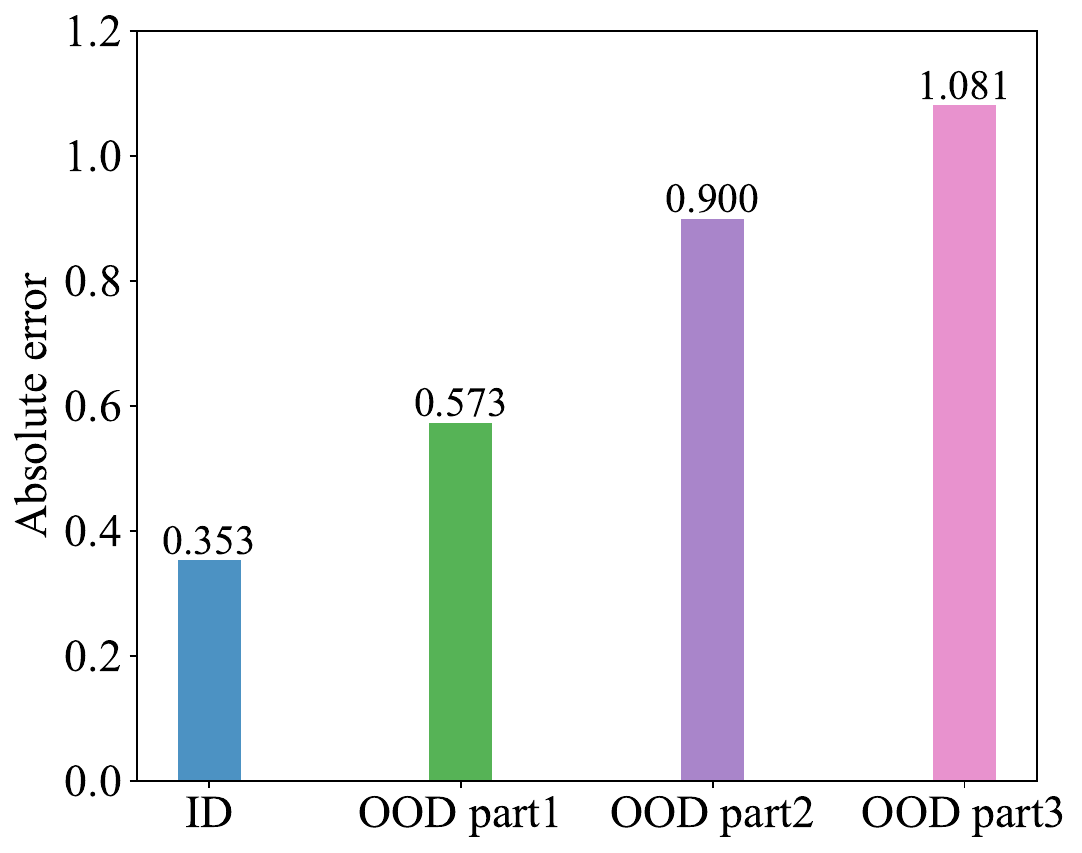}}\hfill
\subfloat[]{\includegraphics[width=0.33\linewidth]{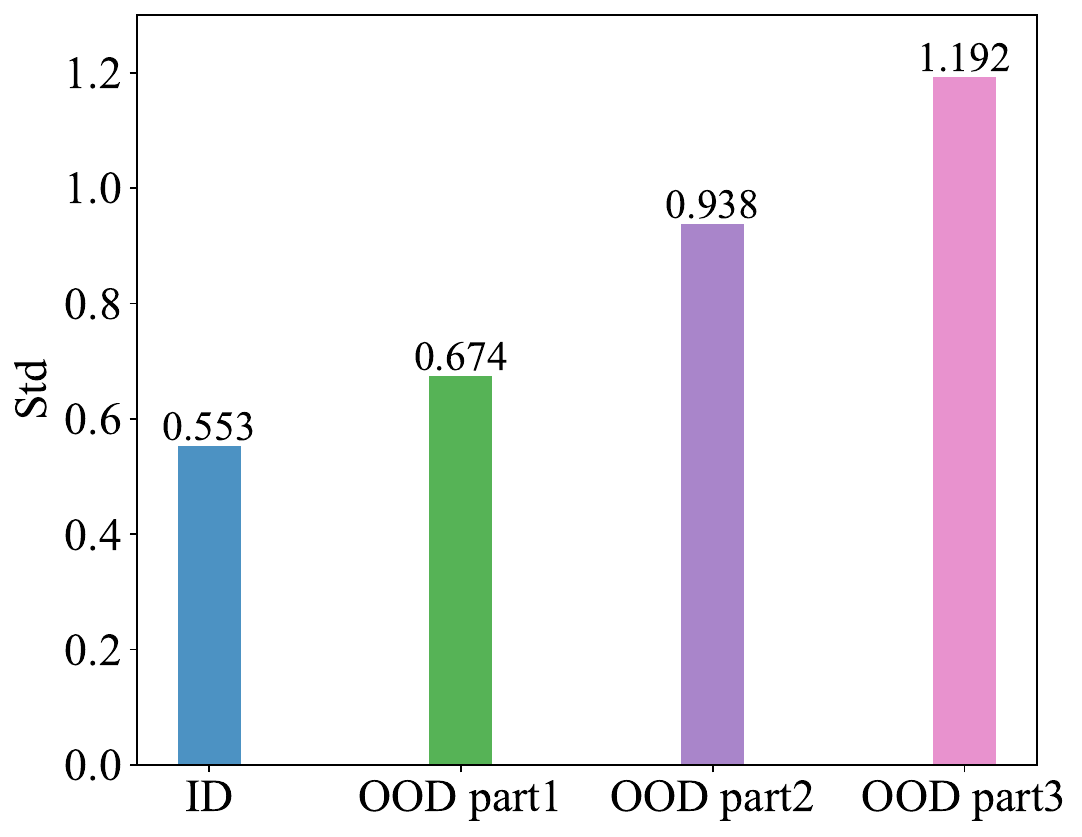}}
\caption{California housing prices regression.  (a) The data partition diagram with the LOF algorithm. Figure (b)-(c) show the absolute error of the predicted mean and standard deviation obtained from IB-UQ, respectively. The absolute error of ID data mean prediction is smaller than OOD data. Further, the absolute error increases significantly for OOD data with decreasing LOF but the uncertainty band covers the absolute error of the mean for both ID and OOD data.}
\label{fig:housing dataset}
\end{figure}

To further evaluate the performance of the proposed IB-UQ method on real data set applications, we next consider the  California housing prices regression problem with data set distributed via the scikit-learn package \cite{pedregosa2011scikit}, which was originally published in \cite{pace1997sparse}. The data set consists of 20640 total samples and 8 features including median household income, median age of housing use, average number of rooms, etc. To obtain the ID and OOD data set, we employ the local outlier factor method (LOF) \cite{breunig2000lof} to divide the whole data set into four blocks based on three  local outlier factor, i.e. -2, -1.5 and -1.2, as shown in the schematic diagram Figure \ref{fig:housing dataset}(a). We regard the corresponding data with LOF score larger than -1.2 as in-distribution data(ID), including 18,090 samples. Data with LOF score smaller than -1.2 is recorded as the OOD data, which is divided into three parts according to the LOF score. Specifically, OOD part1 contains 2192 samples satisfying -1.5<score<-1.2, OOD part2 contains 305 samples satisfying -2<score<-1.5, OOD part3 contains 53 samples satisfying, respectively. We randomly choose $2/3$ of the ID data set for training and the rest samples are used for testing. 

The absolute error of the predict mean and uncertainty estimates for both ID test data and OOD test data are reported in Figure \ref{fig:housing dataset}(b)-(c). We can see from the plots that the absolute error of ID data mean prediction is smaller than those of OOD test data. Further, The absolute error increases significantly for OOD data with decreasing LOF and the uncertainty covers the absolute error of the mean for both ID and OOD data.

%%%%%%%%%%%%%%%%%%%%%%%%%%%%%%%%%%%%%%%%%%%%%%%%%%%%%%%%%%%%%%%%%%%%%%%%%%%%%%%%%%%%%%%%%%%%%%%%%%%%%%%%%%%%%%%%%%

\subsection*{\bf\large{Climate Model}}\label{climate_model}
Finally, we consider a large-scale climate modeling task, which aims to learn an operator that maps the surface air temperature field over the Earth to the surface air pressure field given real historical weather station data \cite{yang2022scalable}. 
This mapping takes the form 
\begin{equation}\label{climate}  
\mathcal{G}: T(x) \rightarrow P(y),
\end{equation}
with $T,P$ represent the temperature and pressure respectively and  $x,y \in [-90, 90] \times [0, 360]$ correspond to latitude and longitude coordinates pairs to specify the location on the surface of the earth. 
\begin{figure}[htp!]
\centering
\subfloat[]{\includegraphics[width=0.99\linewidth]{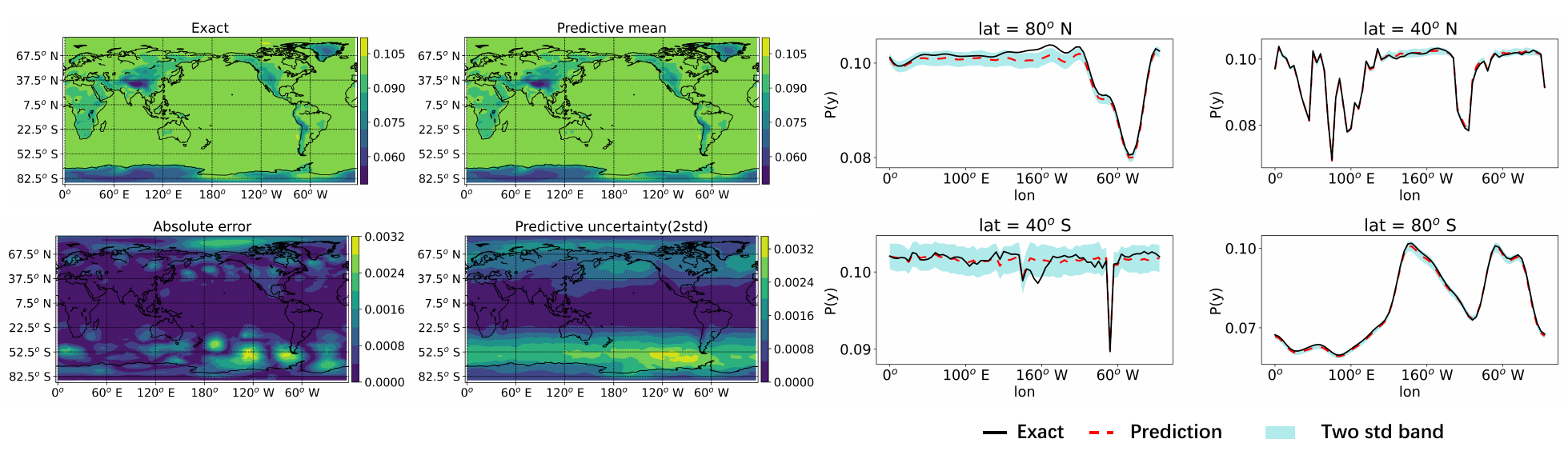}}\\
\subfloat[]{\includegraphics[width=0.99\linewidth]{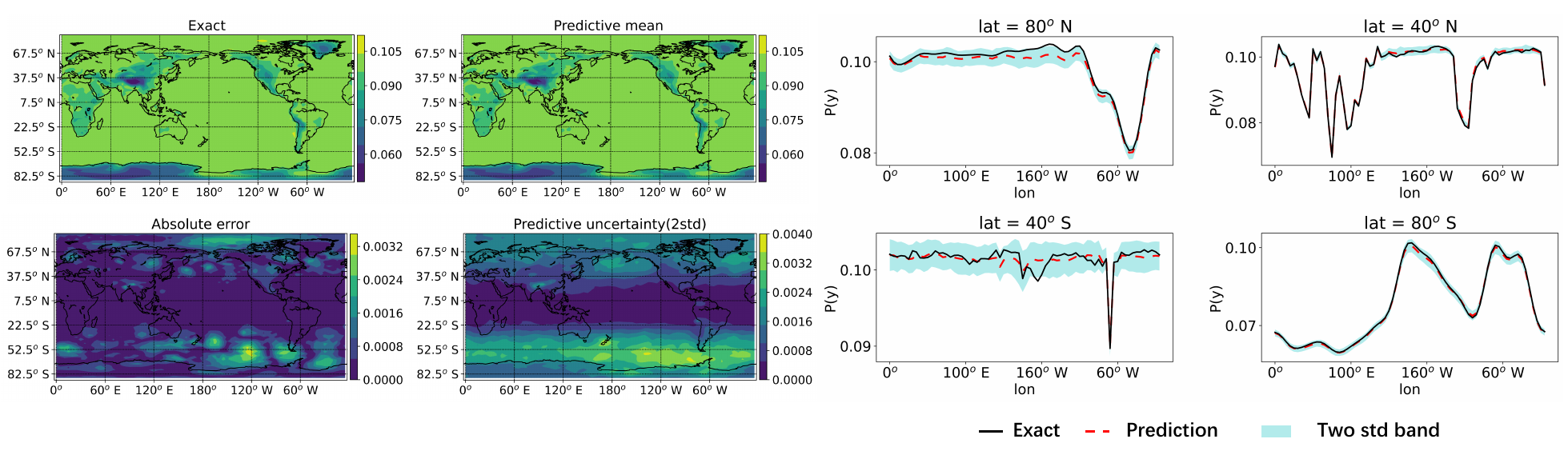}}\\
\subfloat[]{\includegraphics[width=0.99\linewidth]{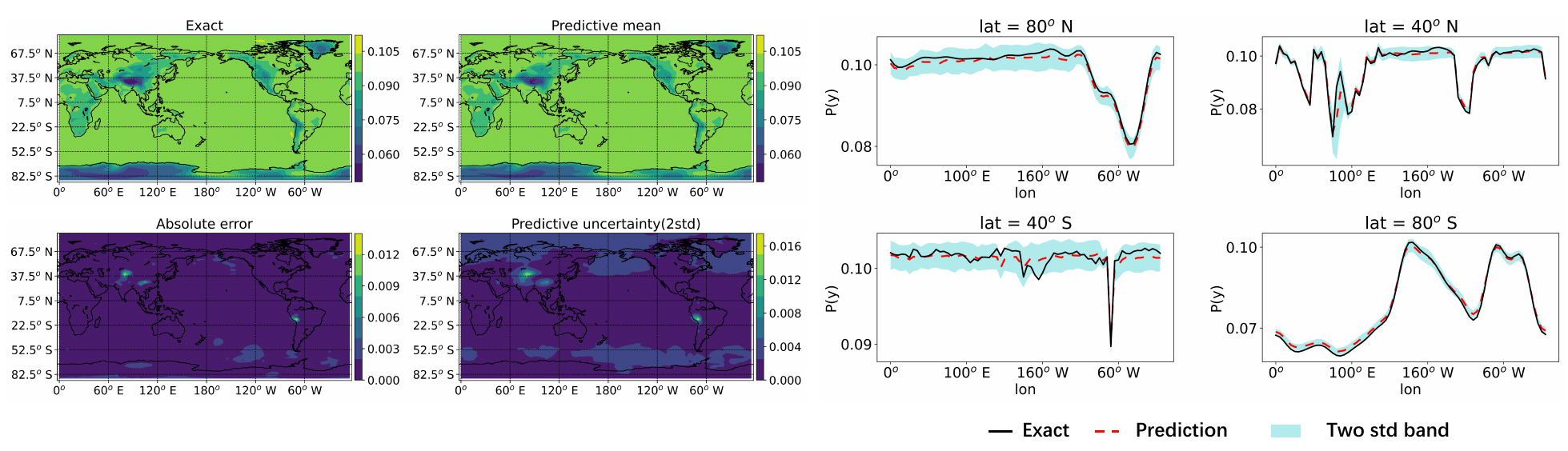}}
\caption{Climate Model. Representative results of the climate model using the IB-UQ algorithm for pressure prediction given surface air temperature. In Figure (a), the forecast results are presented without any removal of training data. The left half of the figure shows the true value of the air pressure, the predicted mean value, the absolute error between the true and predicted values, and twice the standard deviation (forecast uncertainty). The right half of the figure displays the pressure forecast for various latitudes. Figure (b) exhibits the same results as Figure (a), but with all data from the second quarter removed from the training data. Finally, Figure (c) presents the results obtained using training data with all data from the second and third quarters removed.} 
\label{fig:climate_ib}
\end{figure}

The training and testing data sets are obtained according to the data sets used in \cite{yang2022scalable}. The training data set is generated by taking measurements from daily data of surface air temperature and pressure  from 2000 to 2005. Concretely, $N=1825$ input/output function pairs are obtained on a mesh with $144\times 72$ grid. We randomly select $72\times 72$ points during the training stage. The testing data is measured in the same way and contains $N=1825$ input/output function pairs  on a mesh with $72\times 72$ grid. For the details of the data generation, please see the github code accompanying reference \cite{yang2022scalable}.

To enhance the performance of IB-UQ DeepONet for this example, we use following Harmonic Feature Expansion established in \cite{lu2022comprehensive}
\begin{equation}
    \zeta(y) = [y, sin(2\pi y), cos(2\pi y), ... ,  sin(2^{H}\pi y), cos(2^{H}\pi y) ]
\end{equation}
as the input of the trunk network (see Eqs.~(\ref{ONet}, \ref{eq:decoder-deeponet})),
where $H$ is the order of the harmonic basis and we set it equal to 5 in our simulation. Meanwhile, to demonstrate that the IB-UQ method can detect OOD data, we consider three cases in our training process: 1) Case 1 uses full training data set with $N=1825$ sample pairs; 2) Case 2 uses training data after removing the data of the second quarter of each year. The training data set has $N=1375$ sample pairs; 3) Case 3 uses data  after removing all data from the second and third quarter of each year. The training data set for this case has $N=925$ sample pairs. In Figure \ref{fig:climate_ib}, we show the exact pressure field, the predictive mean pressure field, and representative predictions for the three different cases. We observe that the model obtain larger error and uncertainty for the prediction of the unseen data. While for the in distribution data, the model can provide more accurate mean prediction and reasonable uncertainty estimates, which means the uncertainty bound can cover the absolute error.

\FloatBarrier

\section*{\bf\LARGE{Methods}}

\subsection*{\bf\large{IB-UQ for DNN function regression}}
Our proposed approach, termed IB-UQ, utilizes the information bottleneck method for function regression and provides scalable uncertainty estimation in function learning. In this section, we will describe the theory of the method and the architecture of the model in detail.

\textbf{Formulation:} Assume $\mathcal{K}  \subset \mathbb{R}^{d_{X}}$ and $Y=f(X)$ is a function defined on $\mathcal{K}$, i.e. $f:\mathcal{K} \rightarrow \mathbb{R}^{d_{Y}}$. Given a set of paired noisy observations $\mathcal{D} = \{x_i, y_i\}_{i=1}^{N}$, where $x_i\in \mathcal{K}$, $y_i=f(x_i)$. Our goal is to construct a deep learning model that can predict the value of $f$ at any new location $x$ as well as provide reliable uncertainty estimation. Following the main idea of information bottleneck developed by Tishby et al. \cite{tishby2015deep},
we encode the input $X$ to a latent representation $Z$ by enforcing the extraction of essential information relevant to the prediction task and incorporating noise according to the uncertainty caused by the limited training data,
and predict the output $Y$ from $Z$ by a stochastic decoder which can characterize the inherent randomness of the function $f$.
Since the dependence between random variables can be quantified by their mutual information,
the IB principle states that the ideal $Z$ can maximize the IB objective $[I(Z ; Y)-\beta I(Z ; X)]$, where $I(\cdot ; \cdot)$ denotes the mutual information and $\beta\in [0,1]$ is a constant that controls the trade off between prediction and information extraction.

In contrast to the deterministic approach of DNN regression, the IB-UQ approach treats both the latent variable $Z$ and the output $Y$ as random variables.
Moreover, IB-UQ can model both mean and variance in the prediction process and thus provide reliable uncertainty estimation.

\textbf{IB objective with data augmentation:}
In conventional IB methods, $I(Z;Y)$ and $I(Z;X)$ are calculated based on the same joint distribution $(x,y,z)$ given by the training set and encoder. However, such an objective often suffers from the out-of-distribution (OOD) problem. In order to effectively explore areas of the input space with low data distribution density and improve the extrapolation quality,
we utilize a general incompressible-flow network (GIN) to generate
$\tilde X\sim \tilde p(x)$, where $\tilde p(x)\propto p(x)^{\frac{1}{\tau}}$ is a more wide distribution than the training data distribution $p(x)$ with $\tau>1$,
and calculate the mutual information $I(\tilde{Z}; \tilde{X})$ instead of $I(Z;X)$ in the IB objective, i.e.,
\begin{equation}\label{eq:obj-aug}
\mathcal{L}_{\mathrm{IB}}=I(Z;Y)-\beta I(\tilde{Z};\tilde{X})
\end{equation}
As analyzed in Supplementary Information \ref{app:GIN-analysis}, such an objective yields that the latent variable is approximately independent of an OOD input. Moreover, the variational approximation of $I(Z;Y)$ involves the conditional density estimation of $Y|Z$, which may lead to severe overfitting when the training data size is extremely small (see e.g., discussions in \cite{dutordoir2018gaussian,rothfuss2019noise}). In our experiments, we deal with this problem by using Mixup perturbation of the empirical mutual information with small noise intensity. All implementation details of the data augmentation is provided in Supplementary Information \ref{app:gin}.

\begin{figure}[!htp]
\centering
\includegraphics[width=0.8\linewidth]{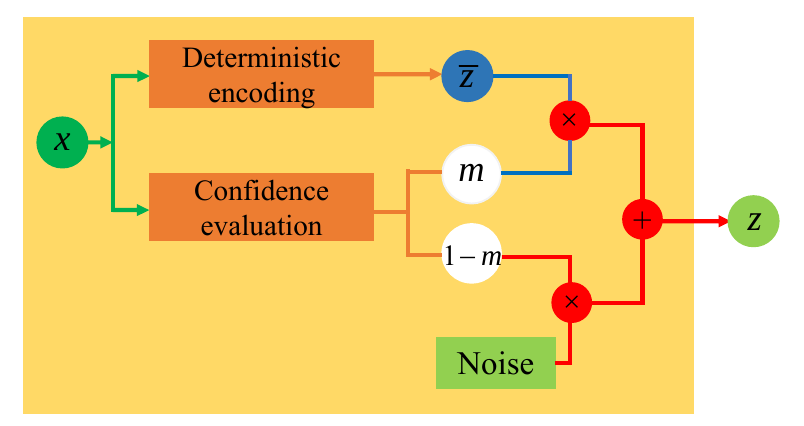}
\caption{Illustration of the encoding process, which involves a deterministic encoding of the input $x$ to $\bar z$ and an evaluation of confidence $m(x)$. The final latent variable $z$ is obtained as a weighted combination of $\bar z$ and noise $z_0$, where the weights are determined by the confidence value $m(x)$.}
\label{fig:Embedding Progress}
\end{figure}

\textbf{Confidence-aware encoder and Gaussian decoder:} Inspired by the stochastic attention mechanism proposed in \cite{mardt2022deep}, the latent representation of the input is given by a confidence-aware encoder
\begin{equation}\label{eq:encoder}
z=\operatorname{diag}(m(x)) \bar{z}(x)+\operatorname{diag}(\mathbf{1}-m(x)) z_{0}
\end{equation}
in IB-UQ as shown in Figure \ref{fig:Embedding Progress}, where $m(x)$, $\bar{z}(x)$ are deep neural networks, each element of $m(x)$ is in [0,1], and $z_{0} \sim \mathcal{N}(0,I)$. Here $\bar z$ is a deterministic and nonlinear feature vector, and  the $i$th element of $m$ represents the confidence of the $i$th feature.
In the case where the encoder is ideally trained, 
$m(x)=\mathbf 1$ implies that $x$ is close to one of training data and the corresponding output $y$ can then be reliably predicted with the deterministic latent variable $z$.
If $m(x) = \mathbf 0$, $x$ is considered as an OOD data and $z$ becomes an uninformative variable $z_{0}$.
Moreover, to appropriately evaluate the confidence for inputs that have very low density in both the training data and GIN model, we set $m(x)$ to be close to zero, which effectively renders the latent variable z uninformative. Further implementation details can be found in the Supplementary Information \ref{app:algorithm}.

The decoder is selected as a $d_y$-dimensional Gaussian model in this paper, i.e., $y$ is conditionally distributed according to $\mathcal{N}(\mu_D(z), \Sigma_D(z))$ for given $z$ in this paper, where $\Sigma_D(z)$ is considered as a diagonal covariance matrix, and $\mu_D(z), \Sigma_D(z)$ are also parameterized as neural networks. It is worth noting that in practical applications, more complex density models such as normalizing flows can be utilized if necessary to capture the distribution of the output more accurately.

\textbf{IB-UQ learning algorithm:}
In general optimizing the IB objective (\ref{eq:obj-aug}) remains difficulty. Here, we utilize the variational approach proposed by  Alemi et al. \cite{alemi2016deep} and derive a tractable lower bound on $\mathcal{L}_{\mathrm{IB}}$ which can be efficiently evaluated:
%maximized using gradient-based methods:
\begin{eqnarray}
\mathcal{L}_{\mathrm{VIB}} & \triangleq & \mathbb{E}_{(x,y,z)\sim p(x,y)\cdot q_{E}(z|x)}\left[\log q_{D}(y|z)\right]-\beta\mathbb{E}_{(\tilde x,\tilde z)\sim \tilde p(\tilde x)\cdot q_{E}(\tilde z|\tilde x)}\left[\log\frac{q_{E}(\tilde z|\tilde x)}{e(\tilde z)}\right]\label{eq:vib_loss}\\
 & \le & \mathcal{L}_{\mathrm{IB}}-\mathcal{H}(Y),\label{eq:vib_bound}
\end{eqnarray}
where $q_E,q_D$ are trainable conditional densities of $Z|X, Y|Z$ given by the encoder and decoder, $e(\tilde z)$ is a parametric model of the marginal distribution of $\tilde Z$ which is employed as normalizing flow (see \ref{app:Realnvpintroduction}) in our IB-UQ model, $\mathcal H(Y)$ denotes the entropy of $Y$, and $p(x,y)$ denotes the joint data distribution. The equality of \eqref{eq:vib_bound} holds if $(q_E,q_D,e)$ are perfectly trained and the maximum value of $\mathcal L_{\mathrm{VIB}}$ is achieved. (See Supplementary Information \ref{app:vib} for derivations.) Therefore, considering that the entropy $\mathcal H(Y)$ of the output is independent of our regression model, we can maximize the variational IB objective $\mathcal L_{\mathrm{VIB}}$ by stochastic gradient descent to find the optimal model parameters. After training, we can predict the mean and uncertainty of $y=f(x)$ for a new $x$ by Monte Carlo simulation of the distribution $\hat p(y|x)=\int q_E(z|x)q_D(y|z)\mathrm dz$.

The schematic of our IB-UQ model is shown in Figure \ref{fig:IBUQ_fun} and implementation details of the proposed algorithm is summarized in Supplementary Information \ref{app:algorithm}.

\subsection*{\bf\large{IB-UQ for neural operator learning}}

In this section, we will establish the information bottleneck method for operator learning. The motivation here is to enhance the robustness of the DeepONet predictions while providing confident uncertainty estimation for out-of-distribution samples.

\textbf{Formulation:} Given $\mathcal{K}_1  \subset \mathbb{R}^{d_{x}}$ and $\mathcal{K}_2 \subset \mathbb{R}^{d_{y}}$, we use  
$C(\mathcal{K}_1; \mathbb{R}^{d_{u}})$ and $C(\mathcal{K}_2; \mathbb{R}^{d_{s}})$ to denote the spaces of continuous input and output functions, $u:\mathcal{K}_1 \rightarrow \mathbb{R}^{d_{u}}$ and $s:\mathcal{K}_2 \rightarrow \mathbb{R}^{d_{s}}$, respectively. $\mathcal{G}:C(\mathcal{K}_1;\mathbb{R}^{d_{u}}) \rightarrow C(\mathcal{K}_2;\mathbb{R}^{d_{s}})$ is a nonlinear operator, $x \in \mathcal{K}_1, y \in \mathcal{K}_2$ are data locations of $u$ and $s$. Suppose we have observations $\{ u^{l}, s^{l} \}_{l=1}^N$ of input and output functions, where  $u^{l} \in C(\mathcal{K}_1; \mathbb{R}^{d_{u}})$, $s^{l} \in C(\mathcal{K}_2; \mathbb{R}^{d_{s}})$, the goal of operator learning is to learn an operator $\mathcal{F}: C(\mathcal{K}_1;\mathbb{R}^{d_{u}}) \rightarrow C(\mathcal{K}_2;\mathbb{R}^{d_{s}})$ satisfying
\begin{equation}\label{sensoreq}
    s^{l} = \mathcal{F}(u^{l}).
\end{equation}

Base on the universal operator approximation theorem developed by Chen \& Chen in \cite{chen1995universal} for single hidden layer neural networks, Lu et.al.~proposed the DeepONet methods for approximating functional and nonlinear operators with deep neural networks in \cite{lu2021learning}. The DeepONet architecture consists of two sub-networks: a trunk network and a branch network. The trunk network takes the data location $y$ in $\mathcal K_2$ as input, while the branch network that takes the sensor measurements $o$ of $u$ as input. Then the outputs of these two sub-networks are combined together by a dot-product operation and finally the output of the DeepONet can be expressed as 
\begin{equation}\label{ONet}
    \mathcal{F}_{\theta}(u)(y) = \sum_{i=1}^{n}b_{i}(o)t_{i}(y)
\end{equation}
where $b$
denotes the $n$ outputs of branch network,
$t$
denotes the $n$ outputs of trunk network, and $\theta$ denotes all trainable parameters of branch and trunk networks, respectively.

Given a data-set consisting of input/output function pairs $\{ u^{l}, s^{l} \}^{N}_{l=1}$, for each pair we assume that we have access to $m$ discrete measurements of the input function at $\{x_j\}_{j=1}^m$ i.e., $o^l = \{u^l(x_j)\}_{j=1}^m$, and $M$ discrete measurements of the output function 
$\{s^l(y_j^l)\}_{j=1}^M$ at query locations $\{y^l_j\}_{j=1}^M$. The original DeepONet model in \cite{lu2021learning} is trained via minimizing the following loss function
\begin{equation}\label{onetl2loss}
    \mathcal{L}(\theta) = \frac{1}{NP} \sum_{l=1}^{N} \sum_{j=1}^{M}\left\vert\mathcal{F}_{\theta}(u^{l})(y^{l}_{j}) - s^{l}(y^{l}_{j}))\right\vert^{2}
    %=\frac{1}{NP} \sum_{l=1}^{N} \sum_{j=1}^{M}\left\vert\sum_{i=1}^nb_i(o^l)t_i(y^{l}_{j}) - s^{l}(y^{l}_{j}))\right\vert^{2}
    .
\end{equation}
Recent findings in \cite{yang2022scalable} and the references therein demonstrate some drawbacks that are inherent to the deterministic training of DeepONet based on (\ref{onetl2loss}). For example, DeepONet can be biased towards approximating
function with large magnitudes. Thus this motivate us to build efficient uncertainty quantification methods for operator learning. 

\textbf{IB-UQ operator learning algorithm}:

We consider to build the information bottleneck in the branch network by assuming that the domain $\mathcal K_2$ of output functions is contained in a low-dimensional space and sufficiently sampled, and construct a stochastic map between the DeepONet input and output functions for UQ by maxmizing the following IB objective
\[
\mathcal L_{\mathrm{IBONet}} = I((Y,Z);s(Y))-\beta I(\tilde Z;\tilde O),
\]
where uppercase letters denote the random variables corresponded to lowercase realizations, $\tilde O$ denotes the synthetic sensor measurements generated by GIN with distribution $\tilde p(o)$, and $\tilde Z$ is the latent variable encoding $\tilde O$.

Now we demonstrate how to compute the variational IB objective of the information bottleneck
based UQ for DeepONet.
We model the encoder as
\begin{equation}\label{eq:encoder-deeponet}
z(o)=\mathrm{diag}(m(o))\bar z(o)+\mathrm{diag}(\mathbf{1}-m(o))z_0
\end{equation}
with $m(o),\bar z(o)$ characterized by neural networks as in \eqref{eq:encoder},
i.e., the conditional distribution of $z$ only depends on input sensor measurements.
The decoder density is also Gaussian with mean $\mu_D(y,z)$ and covariance $\Sigma_D(y,z)=\mathrm{diag}(\sigma_D(y,z))^2$ given by a deterministic DeepONet in the form of
\begin{equation}\label{eq:decoder-deeponet}
\left(\begin{array}{c}
\mu_D(y,z)\\
\log\sigma_D(y,z)
\end{array}\right)=\sum_{i=1}^{n}b_{i}(z)t_{i}(y).
\end{equation}
The encoder and decoder lead to the following variational lower bound of the IB objective
\begin{eqnarray*}
\mathcal{L}_{\mathrm{VIBONet}} & \triangleq & \mathbb{E}_{(o,y,z,s(y))\sim p(o,y,s(y))\cdot q_{E}(z|o)}\left[\log q_{D}(s(y)|y,z)\right]-\beta\mathbb{E}_{(\tilde o,\tilde z)\sim\tilde{p}(\tilde o)\cdot q_{E}(\tilde z|\tilde o)}\left[\log\frac{q_{E}(\tilde z|\tilde o)}{e(\tilde z)}\right],\\
 & \le & \mathcal L_{\mathrm{IBONet}}-\mathcal H(s(Y))
\end{eqnarray*}
with $q_E,q_D$ defined by (\ref{eq:encoder-deeponet}, \ref{eq:decoder-deeponet}).

The schematic of our IB-UQ for DeepONet model is shown in Figure \ref{fig:IBUQ_operator}, and all the derivation and implementation details are given in Supplementary Information \ref{app:operator-algorithm}.

\begin{figure}[htbp]
	\centering
	\includegraphics[width=0.7\linewidth]{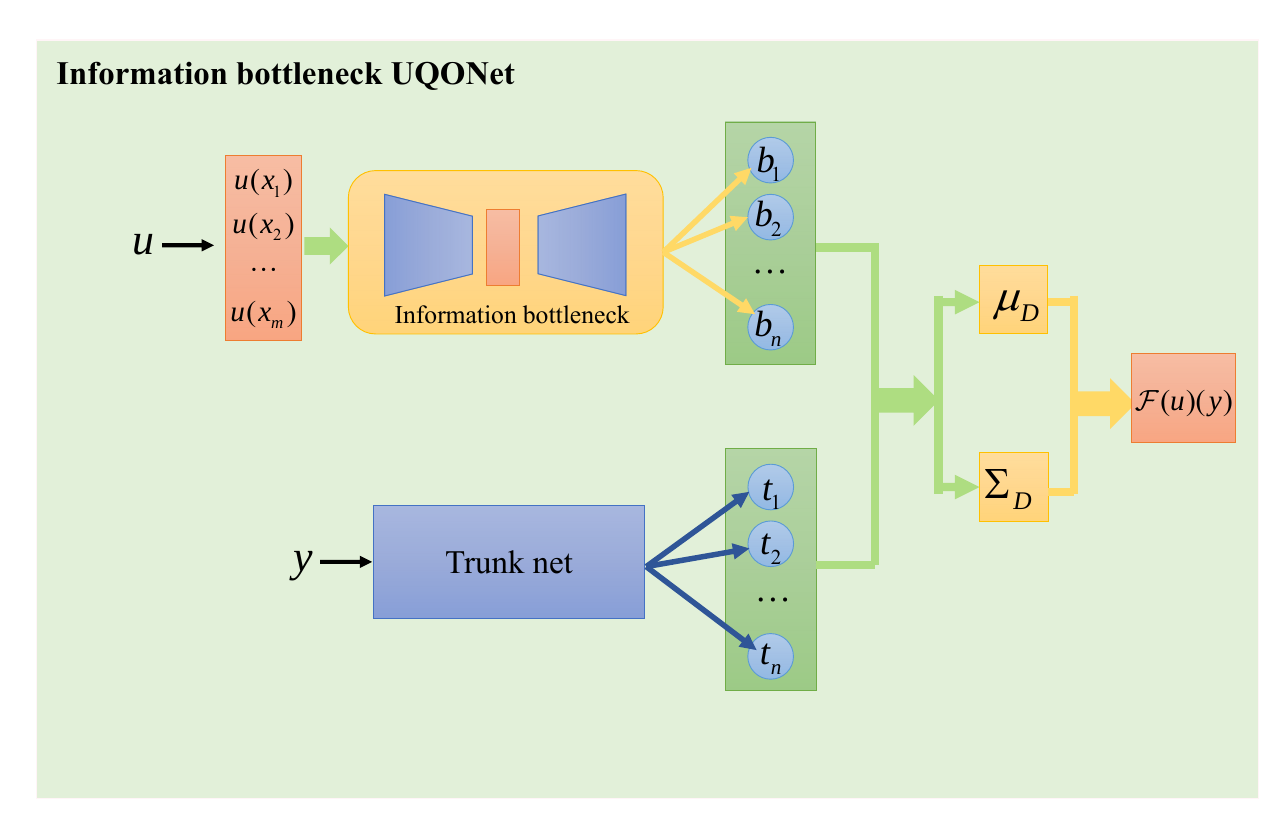}
	\caption{IB-UQ for Operator regression. Schematic of information bottleneck based uncertainty quantification for operator approximation with DeepONet.The Information bottleneck is in the branch of DeepONet. }
	\label{fig:IBUQ_operator}
\end{figure}

\begin{remark}
In this section, we only consider the uncertainty in estimate of $\mathcal F(u)(y)$ at a single location $y$. The UQ for joint prediction values at multiple locations is beyond the scope of this paper and will be investigated in future.
\end{remark}

\begin{remark}
The truck net can also have bottleneck while we choose to use the information bottleneck in the branch network here.
\end{remark}

\section*{\bf\LARGE{Conclusions}}

We have proposed a novel and efficient uncertainty quantification approach based on information bottleneck for function regression and operator learning, coined as IB-UQ. The performance of the proposed IB-UQ methods is tested via several examples, including nonlinear discontinuous function regression, the California housing prices data set and learning the operators of nonlinear partial differential equations. The simulation results demonstrate that the proposed methods can provide robust predictions for function/operator learning as well as reliable uncertainty estimates for noisy limited data and outliers, which is crucial for decision-making in real applications.  

This work is the first attempt to employ information bottleneck for uncertainty quantification in scientific machine learning. In terms of future research, there are several open questions that need further investigation, e.g., extension of the proposed method to high-dimensional problems will be studied in the future work. In addition, while herein we consider only uncertainty quantification for data-driven scientific machine learning, but in practice, we may have partial physics with sparse observations. Thus, IB-UQ for physics informed machine learning need to be investigated. Finally, the capabilities of the proposed IB-UQ methods for active learning for UQ in extreme events and multiscale systems should also be tested and further developed.

\iffalse
\begin{remark}
If we hope to optimize sensor locations of $s$, we can do:
\begin{enumerate}
\item Try different choices of sensor locations.
\item Approximate the entropy of $G$ for all choices (e.g., by normalizing flow), and choose the bottleneck with the best loss.
\end{enumerate}
\end{remark}
\fi

\section*{\bf\large{Acknowledgement}}

We would like to thank Paris Perdikaris (UPenn) for providing the cleaned climate modeling dataset. The first author is supported by the NSF of China (under grant numbers 12071301, 92270115) and the Shanghai Municipal Science and Technology Commission (No. 20JC1412500). The second author is supported by the NSF of China (under grant number 12171367) and the Shanghai Municipal Science and Technology Commission (under grant numbers 20JC1413500 and 2021SHZDZX0100). The last author is supported by the
National Key R\&D Program of China (2020YFA0712000), the NSF of China (under grant numbers 12288201), and youth innovation promotion association (CAS).

\bibliography{Reference} %%#Reference

\begin{thebibliography}{57}
\expandafter\ifx\csname natexlab\endcsname\relax\def\natexlab#1{#1}\fi
\providecommand{\bibinfo}[2]{#2}
\ifx\xfnm\relax \def\xfnm[#1]{\unskip,\space#1}\fi
%Type = Article
\bibitem[{Karniadakis et~al.(2021)Karniadakis, Kevrekidis, Lu, Perdikaris,
  Wang, and Yang}]{karniadakis2021physics}
\bibinfo{author}{G.~E. Karniadakis}, \bibinfo{author}{I.~G. Kevrekidis},
  \bibinfo{author}{L.~Lu}, \bibinfo{author}{P.~Perdikaris},
  \bibinfo{author}{S.~Wang}, \bibinfo{author}{L.~Yang},
\newblock \bibinfo{title}{Physics-informed machine learning},
\newblock \bibinfo{journal}{Nature Reviews Physics} \bibinfo{volume}{3}
  (\bibinfo{year}{2021}) \bibinfo{pages}{422--440}.
%Type = Article
\bibitem[{Willard et~al.(2022)Willard, Jia, Xu, Steinbach, and
  Kumar}]{willard2021integrating}
\bibinfo{author}{J.~Willard}, \bibinfo{author}{X.~Jia},
  \bibinfo{author}{S.~Xu}, \bibinfo{author}{M.~Steinbach},
  \bibinfo{author}{V.~Kumar},
\newblock \bibinfo{title}{Integrating scientific knowledge with machine
  learning for engineering and environmental systems},
\newblock \bibinfo{journal}{ACM Computing Surveys} \bibinfo{volume}{55}
  (\bibinfo{year}{2022}) \bibinfo{pages}{1--37}.
%Type = Article
\bibitem[{Chen and Chen(1993)}]{chen1993approximations}
\bibinfo{author}{T.~Chen}, \bibinfo{author}{H.~Chen},
\newblock \bibinfo{title}{Approximations of continuous functionals by neural
  networks with application to dynamic systems},
\newblock \bibinfo{journal}{IEEE Transactions on Neural networks}
  \bibinfo{volume}{4} (\bibinfo{year}{1993}) \bibinfo{pages}{910--918}.
%Type = Article
\bibitem[{Lagaris et~al.(1998)Lagaris, Likas, and Fotiadis}]{Lagaris1997}
\bibinfo{author}{I.~E. Lagaris}, \bibinfo{author}{A.~C. Likas},
  \bibinfo{author}{D.~I. Fotiadis},
\newblock \bibinfo{title}{Artificial neural networks for solving ordinary and
  partial differential equations},
\newblock \bibinfo{journal}{IEEE Transactions on Neural Networks}
  \bibinfo{volume}{9} (\bibinfo{year}{1998}) \bibinfo{pages}{987--1000}.
%Type = Article
\bibitem[{Khoo et~al.(2021)Khoo, Lu, and Ying}]{khoo2021solving}
\bibinfo{author}{Y.~Khoo}, \bibinfo{author}{J.~Lu}, \bibinfo{author}{L.~Ying},
\newblock \bibinfo{title}{Solving parametric pde problems with artificial
  neural networks},
\newblock \bibinfo{journal}{European Journal of Applied Mathematics}
  \bibinfo{volume}{32} (\bibinfo{year}{2021}) \bibinfo{pages}{421--435}.
%Type = Article
\bibitem[{{Raissi} et~al.(2017){Raissi}, {Perdikaris}, and
  {Karniadakis}}]{MaziarParisGK17_1}
\bibinfo{author}{M.~{Raissi}}, \bibinfo{author}{P.~{Perdikaris}},
  \bibinfo{author}{G.~E. {Karniadakis}},
\newblock \bibinfo{title}{{Physics informed deep learning (part I): Data-driven
  solutions of nonlinear partial differential equations}},
\newblock \bibinfo{journal}{arXiv preprint}  (\bibinfo{year}{2017})
  \bibinfo{pages}{arXiv:1711.10561}.
%Type = Article
\bibitem[{Liao and Ming(2021)}]{LiaoMing2021}
\bibinfo{author}{Y.~Liao}, \bibinfo{author}{P.~Ming},
\newblock \bibinfo{title}{Deep {N}itsche method: Deep {R}itz method with
  essential boundary conditions},
\newblock \bibinfo{journal}{Commun. Comput. Phys, 29:1365-1384}
  (\bibinfo{year}{2021}).
%Type = Article
\bibitem[{Guo et~al.(2022{\natexlab{a}})Guo, Wu, and Zhou}]{guo2022normalizing}
\bibinfo{author}{L.~Guo}, \bibinfo{author}{H.~Wu}, \bibinfo{author}{T.~Zhou},
\newblock \bibinfo{title}{Normalizing field flows: Solving forward and inverse
  stochastic differential equations using physics-informed flow models},
\newblock \bibinfo{journal}{Journal of Computational Physics}
  \bibinfo{volume}{461} (\bibinfo{year}{2022}{\natexlab{a}})
  \bibinfo{pages}{111202}.
%Type = Article
\bibitem[{Guo et~al.(2022{\natexlab{b}})Guo, Wu, Yu, and Zhou}]{guo2022monte}
\bibinfo{author}{L.~Guo}, \bibinfo{author}{H.~Wu}, \bibinfo{author}{X.~Yu},
  \bibinfo{author}{T.~Zhou},
\newblock \bibinfo{title}{Monte carlo pinns: deep learning approach for forward
  and inverse problems involving high dimensional fractional partial
  differential equations},
\newblock \bibinfo{journal}{arXiv preprint arXiv:2203.08501}
  (\bibinfo{year}{2022}{\natexlab{b}}).
%Type = Article
\bibitem[{Huang et~al.(2022)Huang, Wang, and Zhou}]{huang2022augmented}
\bibinfo{author}{J.~Huang}, \bibinfo{author}{H.~Wang},
  \bibinfo{author}{T.~Zhou},
\newblock \bibinfo{title}{An augmented lagrangian deep learning method for
  variational problems with essential boundary conditions},
\newblock \bibinfo{journal}{Communications in Computational Physics}
  \bibinfo{volume}{31} (\bibinfo{year}{2022}) \bibinfo{pages}{966--986}.
%Type = Article
\bibitem[{Gao et~al.(2022)Gao, Yan, and Zhou}]{gao2022failure}
\bibinfo{author}{Z.~Gao}, \bibinfo{author}{L.~Yan}, \bibinfo{author}{T.~Zhou},
\newblock \bibinfo{title}{Failure-informed adaptive sampling for pinns},
\newblock \bibinfo{journal}{arXiv preprint arXiv:2210.00279}
  (\bibinfo{year}{2022}).
%Type = Article
\bibitem[{E and Yu(2018)}]{EYu2018}
\bibinfo{author}{W.~E}, \bibinfo{author}{B.~Yu},
\newblock \bibinfo{title}{The deep {R}itz method: a deep learning-based
  numerical algorithm for solving variational problems},
\newblock \bibinfo{journal}{Commun. Math. Stat., 6:1-12}
  (\bibinfo{year}{2018}).
%Type = Article
\bibitem[{Zang et~al.(2020)Zang, Bao, Ye, and Zhou}]{ZangBaoYeZhou2020}
\bibinfo{author}{Y.~Zang}, \bibinfo{author}{G.~Bao}, \bibinfo{author}{X.~Ye},
  \bibinfo{author}{H.~Zhou},
\newblock \bibinfo{title}{Weak adversarial networks for high-dimensional
  partial differential equations},
\newblock \bibinfo{journal}{J. Comput. Phys., 411:109409}
  (\bibinfo{year}{2020}).
%Type = Article
\bibitem[{Brunton et~al.(2016)Brunton, Proctor, and
  Kutz}]{brunton2016discovering}
\bibinfo{author}{S.~L. Brunton}, \bibinfo{author}{J.~L. Proctor},
  \bibinfo{author}{J.~N. Kutz},
\newblock \bibinfo{title}{Discovering governing equations from data by sparse
  identification of nonlinear dynamical systems},
\newblock \bibinfo{journal}{Proceedings of the national academy of sciences}
  \bibinfo{volume}{113} (\bibinfo{year}{2016}) \bibinfo{pages}{3932--3937}.
%Type = Inproceedings
\bibitem[{Long et~al.(????)Long, Lu, Ma, and Dong}]{long2018pde}
\bibinfo{author}{Z.~Long}, \bibinfo{author}{Y.~Lu}, \bibinfo{author}{X.~Ma},
  \bibinfo{author}{B.~Dong},
\newblock \bibinfo{title}{Pde-net: Learning pdes from data},
\newblock in: \bibinfo{booktitle}{International conference on machine learning
  (2018) 3208--3216}, \bibinfo{organization}{PMLR}.
%Type = Article
\bibitem[{Lu et~al.(2021)Lu, Jin, Pang, Zhang, and
  Karniadakis}]{lu2021learning}
\bibinfo{author}{L.~Lu}, \bibinfo{author}{P.~Jin}, \bibinfo{author}{G.~Pang},
  \bibinfo{author}{Z.~Zhang}, \bibinfo{author}{G.~E. Karniadakis},
\newblock \bibinfo{title}{Learning nonlinear operators via deeponet based on
  the universal approximation theorem of operators},
\newblock \bibinfo{journal}{Nature Machine Intelligence} \bibinfo{volume}{3}
  (\bibinfo{year}{2021}) \bibinfo{pages}{218--229}.
%Type = Article
\bibitem[{Chen and Chen(1995)}]{chen1995universal}
\bibinfo{author}{T.~Chen}, \bibinfo{author}{H.~Chen},
\newblock \bibinfo{title}{Universal approximation to nonlinear operators by
  neural networks with arbitrary activation functions and its application to
  dynamical systems},
\newblock \bibinfo{journal}{IEEE Transactions on Neural Networks}
  \bibinfo{volume}{6} (\bibinfo{year}{1995}) \bibinfo{pages}{911--917}.
%Type = Article
\bibitem[{Lin et~al.(2021)Lin, Li, Lu, Cai, Maxey, and
  Karniadakis}]{lin2021operator}
\bibinfo{author}{C.~Lin}, \bibinfo{author}{Z.~Li}, \bibinfo{author}{L.~Lu},
  \bibinfo{author}{S.~Cai}, \bibinfo{author}{M.~Maxey}, \bibinfo{author}{G.~E.
  Karniadakis},
\newblock \bibinfo{title}{Operator learning for predicting multiscale bubble
  growth dynamics},
\newblock \bibinfo{journal}{The Journal of Chemical Physics}
  \bibinfo{volume}{154} (\bibinfo{year}{2021}) \bibinfo{pages}{104118}.
%Type = Article
\bibitem[{Mao et~al.(2021)Mao, Lu, Marxen, Zaki, and
  Karniadakis}]{mao2021deepm}
\bibinfo{author}{Z.~Mao}, \bibinfo{author}{L.~Lu}, \bibinfo{author}{O.~Marxen},
  \bibinfo{author}{T.~A. Zaki}, \bibinfo{author}{G.~E. Karniadakis},
\newblock \bibinfo{title}{Deepm\&mnet for hypersonics: Predicting the coupled
  flow and finite-rate chemistry behind a normal shock using neural-network
  approximation of operators},
\newblock \bibinfo{journal}{Journal of computational physics}
  \bibinfo{volume}{447} (\bibinfo{year}{2021}) \bibinfo{pages}{110698}.
%Type = Article
\bibitem[{Wang et~al.(2021)Wang, Wang, and Perdikaris}]{wang2021learning}
\bibinfo{author}{S.~Wang}, \bibinfo{author}{H.~Wang},
  \bibinfo{author}{P.~Perdikaris},
\newblock \bibinfo{title}{Learning the solution operator of parametric partial
  differential equations with physics-informed deeponets},
\newblock \bibinfo{journal}{Science advances} \bibinfo{volume}{7}
  (\bibinfo{year}{2021}) \bibinfo{pages}{eabi8605}.
%Type = Inproceedings
\bibitem[{Li et~al.(2020)Li, Kovachki, Azizzadenesheli, Bhattacharya, Stuart,
  Anandkumar et~al.}]{li2020fourier}
\bibinfo{author}{Z.~Li}, \bibinfo{author}{N.~B. Kovachki},
  \bibinfo{author}{K.~Azizzadenesheli}, \bibinfo{author}{K.~Bhattacharya},
  \bibinfo{author}{A.~Stuart}, \bibinfo{author}{A.~Anandkumar}, et~al.,
\newblock \bibinfo{title}{Fourier neural operator for parametric partial
  differential equations},
\newblock in: \bibinfo{booktitle}{International Conference on Learning
  Representations (2020)}.
%Type = Article
\bibitem[{Kissas et~al.(2022)Kissas, Seidman, Guilhoto, Preciado, Pappas, and
  Perdikaris}]{kissas2022learning}
\bibinfo{author}{G.~Kissas}, \bibinfo{author}{J.~H. Seidman},
  \bibinfo{author}{L.~F. Guilhoto}, \bibinfo{author}{V.~M. Preciado},
  \bibinfo{author}{G.~J. Pappas}, \bibinfo{author}{P.~Perdikaris},
\newblock \bibinfo{title}{Learning operators with coupled attention},
\newblock \bibinfo{journal}{Journal of Machine Learning Research}
  \bibinfo{volume}{23} (\bibinfo{year}{2022}) \bibinfo{pages}{1--63}.
%Type = Article
\bibitem[{Psaros et~al.(2022)Psaros, Meng, Zou, Guo, and
  Karniadakis}]{psaros2022uncertainty}
\bibinfo{author}{A.~F. Psaros}, \bibinfo{author}{X.~Meng},
  \bibinfo{author}{Z.~Zou}, \bibinfo{author}{L.~Guo}, \bibinfo{author}{G.~E.
  Karniadakis},
\newblock \bibinfo{title}{Uncertainty quantification in scientific machine
  learning: Methods, metrics, and comparisons},
\newblock \bibinfo{journal}{arXiv preprint arXiv:2201.07766}
  (\bibinfo{year}{2022}).
%Type = Article
\bibitem[{MacKay(1995)}]{mackay1995bayesian}
\bibinfo{author}{D.~J. MacKay},
\newblock \bibinfo{title}{Bayesian neural networks and density networks},
\newblock \bibinfo{journal}{Nuclear Instruments and Methods in Physics Research
  Section A: Accelerators, Spectrometers, Detectors and Associated Equipment}
  \bibinfo{volume}{354} (\bibinfo{year}{1995}) \bibinfo{pages}{73--80}.
%Type = Book
\bibitem[{Neal(2012)}]{neal2012bayesian}
\bibinfo{author}{R.~M. Neal}, \bibinfo{title}{Bayesian learning for neural
  networks}, volume \bibinfo{volume}{118}, \bibinfo{publisher}{Springer Science
  \& Business Media}, \bibinfo{year}{2012}.
%Type = Inproceedings
\bibitem[{Gal and Ghahramani(????)}]{gal2016dropout}
\bibinfo{author}{Y.~Gal}, \bibinfo{author}{Z.~Ghahramani},
\newblock \bibinfo{title}{Dropout as a bayesian approximation: Representing
  model uncertainty in deep learning},
\newblock in: \bibinfo{booktitle}{international conference on machine learning
  (2016) 1050--1059}, \bibinfo{organization}{PMLR}.
%Type = Article
\bibitem[{Lakshminarayanan et~al.(2017)Lakshminarayanan, Pritzel, and
  Blundell}]{lakshminarayanan2017simple}
\bibinfo{author}{B.~Lakshminarayanan}, \bibinfo{author}{A.~Pritzel},
  \bibinfo{author}{C.~Blundell},
\newblock \bibinfo{title}{Simple and scalable predictive uncertainty estimation
  using deep ensembles},
\newblock \bibinfo{journal}{Advances in neural information processing systems}
  \bibinfo{volume}{30} (\bibinfo{year}{2017}).
%Type = Article
\bibitem[{Malinin and Gales(2018)}]{malinin2018predictive}
\bibinfo{author}{A.~Malinin}, \bibinfo{author}{M.~Gales},
\newblock \bibinfo{title}{Predictive uncertainty estimation via prior
  networks},
\newblock \bibinfo{journal}{Advances in neural information processing systems}
  \bibinfo{volume}{31} (\bibinfo{year}{2018}).
%Type = Article
\bibitem[{Fort et~al.(2019)Fort, Hu, and Lakshminarayanan}]{fort2019deep}
\bibinfo{author}{S.~Fort}, \bibinfo{author}{H.~Hu},
  \bibinfo{author}{B.~Lakshminarayanan},
\newblock \bibinfo{title}{Deep ensembles: A loss landscape perspective},
\newblock \bibinfo{journal}{arXiv preprint arXiv:1912.02757}
  (\bibinfo{year}{2019}).
%Type = Article
\bibitem[{Sensoy et~al.(2018)Sensoy, Kaplan, and
  Kandemir}]{sensoy2018evidential}
\bibinfo{author}{M.~Sensoy}, \bibinfo{author}{L.~Kaplan},
  \bibinfo{author}{M.~Kandemir},
\newblock \bibinfo{title}{Evidential deep learning to quantify classification
  uncertainty},
\newblock \bibinfo{journal}{Advances in neural information processing systems}
  \bibinfo{volume}{31} (\bibinfo{year}{2018}).
%Type = Article
\bibitem[{Amini et~al.(2020)Amini, Schwarting, Soleimany, and
  Rus}]{amini2020deep}
\bibinfo{author}{A.~Amini}, \bibinfo{author}{W.~Schwarting},
  \bibinfo{author}{A.~Soleimany}, \bibinfo{author}{D.~Rus},
\newblock \bibinfo{title}{Deep evidential regression},
\newblock \bibinfo{journal}{Advances in Neural Information Processing Systems}
  \bibinfo{volume}{33} (\bibinfo{year}{2020}) \bibinfo{pages}{14927--14937}.
%Type = Article
\bibitem[{Zou et~al.(2022)Zou, Meng, Psaros, and Karniadakis}]{zou2022neuraluq}
\bibinfo{author}{Z.~Zou}, \bibinfo{author}{X.~Meng}, \bibinfo{author}{A.~F.
  Psaros}, \bibinfo{author}{G.~E. Karniadakis},
\newblock \bibinfo{title}{Neuraluq: A comprehensive library for uncertainty
  quantification in neural differential equations and operators},
\newblock \bibinfo{journal}{arXiv preprint arXiv:2208.11866}
  (\bibinfo{year}{2022}).
%Type = Article
\bibitem[{Moya et~al.(2023)Moya, Zhang, Lin, and Yue}]{moya2023deeponet}
\bibinfo{author}{C.~Moya}, \bibinfo{author}{S.~Zhang},
  \bibinfo{author}{G.~Lin}, \bibinfo{author}{M.~Yue},
\newblock \bibinfo{title}{Deeponet-grid-uq: A trustworthy deep operator
  framework for predicting the power grid’s post-fault trajectories},
\newblock \bibinfo{journal}{Neurocomputing} \bibinfo{volume}{535}
  (\bibinfo{year}{2023}) \bibinfo{pages}{166--182}.
%Type = Article
\bibitem[{Lin et~al.(2021)Lin, Moya, and Zhang}]{lin2021accelerated}
\bibinfo{author}{G.~Lin}, \bibinfo{author}{C.~Moya},
  \bibinfo{author}{Z.~Zhang},
\newblock \bibinfo{title}{Accelerated replica exchange stochastic gradient
  langevin diffusion enhanced bayesian deeponet for solving noisy parametric
  pdes},
\newblock \bibinfo{journal}{arXiv preprint arXiv:2111.02484}
  (\bibinfo{year}{2021}).
%Type = Article
\bibitem[{Yang et~al.(2022)Yang, Kissas, and Perdikaris}]{yang2022scalable}
\bibinfo{author}{Y.~Yang}, \bibinfo{author}{G.~Kissas},
  \bibinfo{author}{P.~Perdikaris},
\newblock \bibinfo{title}{Scalable uncertainty quantification for deep operator
  networks using randomized priors},
\newblock \bibinfo{journal}{Computer Methods in Applied Mechanics and
  Engineering} \bibinfo{volume}{399} (\bibinfo{year}{2022})
  \bibinfo{pages}{115399}.
%Type = Article
\bibitem[{Zhu et~al.(2022)Zhu, Zhang, Jiao, Karniadakis, and
  Lu}]{zhu2022reliable}
\bibinfo{author}{M.~Zhu}, \bibinfo{author}{H.~Zhang},
  \bibinfo{author}{A.~Jiao}, \bibinfo{author}{G.~E. Karniadakis},
  \bibinfo{author}{L.~Lu},
\newblock \bibinfo{title}{Reliable extrapolation of deep neural operators
  informed by physics or sparse observations},
\newblock \bibinfo{journal}{arXiv preprint arXiv:2212.06347}
  (\bibinfo{year}{2022}).
%Type = Inproceedings
\bibitem[{Tishby(1999)}]{tishby1999information}
\bibinfo{author}{N.~Tishby},
\newblock \bibinfo{title}{The information bottleneck method},
\newblock in: \bibinfo{booktitle}{Proc. 37th Annual Allerton Conference on
  Communications, Control and Computing (1999)}.
%Type = Inproceedings
\bibitem[{Tishby and Zaslavsky(????)}]{tishby2015deep}
\bibinfo{author}{N.~Tishby}, \bibinfo{author}{N.~Zaslavsky},
\newblock \bibinfo{title}{Deep learning and the information bottleneck
  principle},
\newblock in: \bibinfo{booktitle}{ieee information theory workshop (itw) (2015)
  1--5}.
%Type = Article
\bibitem[{Alemi et~al.(2016)Alemi, Fischer, Dillon, and Murphy}]{alemi2016deep}
\bibinfo{author}{A.~A. Alemi}, \bibinfo{author}{I.~Fischer},
  \bibinfo{author}{J.~V. Dillon}, \bibinfo{author}{K.~Murphy},
\newblock \bibinfo{title}{Deep variational information bottleneck},
\newblock \bibinfo{journal}{arXiv preprint arXiv:1612.00410}
  (\bibinfo{year}{2016}).
%Type = Article
\bibitem[{Alemi et~al.(2018)Alemi, Fischer, and Dillon}]{alemi2018uncertainty}
\bibinfo{author}{A.~A. Alemi}, \bibinfo{author}{I.~Fischer},
  \bibinfo{author}{J.~V. Dillon},
\newblock \bibinfo{title}{Uncertainty in the variational information
  bottleneck},
\newblock \bibinfo{journal}{arXiv preprint arXiv:1807.00906}
  (\bibinfo{year}{2018}).
%Type = Article
\bibitem[{Kolchinsky et~al.(2019)Kolchinsky, Tracey, and
  Wolpert}]{kolchinsky2019nonlinear}
\bibinfo{author}{A.~Kolchinsky}, \bibinfo{author}{B.~D. Tracey},
  \bibinfo{author}{D.~H. Wolpert},
\newblock \bibinfo{title}{Nonlinear information bottleneck},
\newblock \bibinfo{journal}{Entropy} \bibinfo{volume}{21}
  (\bibinfo{year}{2019}) \bibinfo{pages}{1181}.
%Type = Article
\bibitem[{Shwartz-Ziv and Tishby(2017)}]{shwartz2017opening}
\bibinfo{author}{R.~Shwartz-Ziv}, \bibinfo{author}{N.~Tishby},
\newblock \bibinfo{title}{Opening the black box of deep neural networks via
  information},
\newblock \bibinfo{journal}{arXiv preprint arXiv:1703.00810}
  (\bibinfo{year}{2017}).
%Type = Article
\bibitem[{Saxe et~al.(2019)Saxe, Bansal, Dapello, Advani, Kolchinsky, Tracey,
  and Cox}]{saxe2019information}
\bibinfo{author}{A.~M. Saxe}, \bibinfo{author}{Y.~Bansal},
  \bibinfo{author}{J.~Dapello}, \bibinfo{author}{M.~Advani},
  \bibinfo{author}{A.~Kolchinsky}, \bibinfo{author}{B.~D. Tracey},
  \bibinfo{author}{D.~D. Cox},
\newblock \bibinfo{title}{On the information bottleneck theory of deep
  learning},
\newblock \bibinfo{journal}{Journal of Statistical Mechanics: Theory and
  Experiment} \bibinfo{volume}{2019} (\bibinfo{year}{2019})
  \bibinfo{pages}{124020}.
%Type = Article
\bibitem[{Pedregosa et~al.(2011)Pedregosa, Varoquaux, Gramfort, Michel,
  Thirion, Grisel, Blondel, Prettenhofer, Weiss, Dubourg
  et~al.}]{pedregosa2011scikit}
\bibinfo{author}{F.~Pedregosa}, \bibinfo{author}{G.~Varoquaux},
  \bibinfo{author}{A.~Gramfort}, \bibinfo{author}{V.~Michel},
  \bibinfo{author}{B.~Thirion}, \bibinfo{author}{O.~Grisel},
  \bibinfo{author}{M.~Blondel}, \bibinfo{author}{P.~Prettenhofer},
  \bibinfo{author}{R.~Weiss}, \bibinfo{author}{V.~Dubourg}, et~al.,
\newblock \bibinfo{title}{Scikit-learn: Machine learning in python},
\newblock \bibinfo{journal}{the Journal of machine Learning research}
  \bibinfo{volume}{12} (\bibinfo{year}{2011}) \bibinfo{pages}{2825--2830}.
%Type = Article
\bibitem[{Pace and Barry(1997)}]{pace1997sparse}
\bibinfo{author}{R.~K. Pace}, \bibinfo{author}{R.~Barry},
\newblock \bibinfo{title}{Sparse spatial autoregressions},
\newblock \bibinfo{journal}{Statistics \& Probability Letters}
  \bibinfo{volume}{33} (\bibinfo{year}{1997}) \bibinfo{pages}{291--297}.
%Type = Inproceedings
\bibitem[{Breunig et~al.(????)Breunig, Kriegel, Ng, and
  Sander}]{breunig2000lof}
\bibinfo{author}{M.~M. Breunig}, \bibinfo{author}{H.-P. Kriegel},
  \bibinfo{author}{R.~T. Ng}, \bibinfo{author}{J.~Sander},
\newblock \bibinfo{title}{Lof: identifying density-based local outliers},
\newblock in: \bibinfo{booktitle}{Proceedings of the 2000 ACM SIGMOD
  international conference on Management of data (2000) 93--104}.
%Type = Article
\bibitem[{Lu et~al.(2022)Lu, Meng, Cai, Mao, Goswami, Zhang, and
  Karniadakis}]{lu2022comprehensive}
\bibinfo{author}{L.~Lu}, \bibinfo{author}{X.~Meng}, \bibinfo{author}{S.~Cai},
  \bibinfo{author}{Z.~Mao}, \bibinfo{author}{S.~Goswami},
  \bibinfo{author}{Z.~Zhang}, \bibinfo{author}{G.~E. Karniadakis},
\newblock \bibinfo{title}{A comprehensive and fair comparison of two neural
  operators (with practical extensions) based on fair data},
\newblock \bibinfo{journal}{Computer Methods in Applied Mechanics and
  Engineering} \bibinfo{volume}{393} (\bibinfo{year}{2022})
  \bibinfo{pages}{114778}.
%Type = Article
\bibitem[{Dutordoir et~al.(2018)Dutordoir, Salimbeni, Hensman, and
  Deisenroth}]{dutordoir2018gaussian}
\bibinfo{author}{V.~Dutordoir}, \bibinfo{author}{H.~Salimbeni},
  \bibinfo{author}{J.~Hensman}, \bibinfo{author}{M.~Deisenroth},
\newblock \bibinfo{title}{Gaussian process conditional density estimation},
\newblock \bibinfo{journal}{Advances in neural information processing systems}
  \bibinfo{volume}{31} (\bibinfo{year}{2018}).
%Type = Article
\bibitem[{Rothfuss et~al.(2019)Rothfuss, Ferreira, Boehm, Walther, Ulrich,
  Asfour, and Krause}]{rothfuss2019noise}
\bibinfo{author}{J.~Rothfuss}, \bibinfo{author}{F.~Ferreira},
  \bibinfo{author}{S.~Boehm}, \bibinfo{author}{S.~Walther},
  \bibinfo{author}{M.~Ulrich}, \bibinfo{author}{T.~Asfour},
  \bibinfo{author}{A.~Krause},
\newblock \bibinfo{title}{Noise regularization for conditional density
  estimation},
\newblock \bibinfo{journal}{arXiv preprint arXiv:1907.08982}
  (\bibinfo{year}{2019}).
%Type = Article
\bibitem[{Mardt et~al.(2022)Mardt, Hempel, Clementi, and
  No{\'e}}]{mardt2022deep}
\bibinfo{author}{A.~Mardt}, \bibinfo{author}{T.~Hempel},
  \bibinfo{author}{C.~Clementi}, \bibinfo{author}{F.~No{\'e}},
\newblock \bibinfo{title}{Deep learning to decompose macromolecules into
  independent markovian domains},
\newblock \bibinfo{journal}{Nature Communications} \bibinfo{volume}{13}
  (\bibinfo{year}{2022}) \bibinfo{pages}{7101}.
%Type = Article
\bibitem[{Dinh et~al.(2016)Dinh, Sohl-Dickstein, and Bengio}]{dinh2016density}
\bibinfo{author}{L.~Dinh}, \bibinfo{author}{J.~Sohl-Dickstein},
  \bibinfo{author}{S.~Bengio},
\newblock \bibinfo{title}{Density estimation using real nvp},
\newblock \bibinfo{journal}{arXiv preprint arXiv:1605.08803}
  (\bibinfo{year}{2016}).
%Type = Article
\bibitem[{Kobyzev et~al.(2020)Kobyzev, Prince, and
  Brubaker}]{kobyzev2020normalizing}
\bibinfo{author}{I.~Kobyzev}, \bibinfo{author}{S.~J. Prince},
  \bibinfo{author}{M.~A. Brubaker},
\newblock \bibinfo{title}{Normalizing flows: An introduction and review of
  current methods},
\newblock \bibinfo{journal}{IEEE transactions on pattern analysis and machine
  intelligence} \bibinfo{volume}{43} (\bibinfo{year}{2020})
  \bibinfo{pages}{3964--3979}.
%Type = Inproceedings
\bibitem[{Sorrenson et~al.(2020)Sorrenson, Rother, and
  K{\"o}the}]{sorrenson2019disentanglement}
\bibinfo{author}{P.~Sorrenson}, \bibinfo{author}{C.~Rother},
  \bibinfo{author}{U.~K{\"o}the},
\newblock \bibinfo{title}{Disentanglement by nonlinear ica with general
  incompressible-flow networks (gin)},
\newblock in: \bibinfo{booktitle}{International Conference on Learning
  Representations (2020)}.
%Type = Article
\bibitem[{Dibak et~al.(2022)Dibak, Klein, Kr{\"a}mer, and
  No{\'e}}]{dibak2022temperature}
\bibinfo{author}{M.~Dibak}, \bibinfo{author}{L.~Klein},
  \bibinfo{author}{A.~Kr{\"a}mer}, \bibinfo{author}{F.~No{\'e}},
\newblock \bibinfo{title}{Temperature steerable flows and boltzmann
  generators},
\newblock \bibinfo{journal}{Physical Review Research} \bibinfo{volume}{4}
  (\bibinfo{year}{2022}) \bibinfo{pages}{L042005}.
%Type = Inproceedings
\bibitem[{Cornish et~al.(????)Cornish, Caterini, Deligiannidis, and
  Doucet}]{cornish2020relaxing}
\bibinfo{author}{R.~Cornish}, \bibinfo{author}{A.~Caterini},
  \bibinfo{author}{G.~Deligiannidis}, \bibinfo{author}{A.~Doucet},
\newblock \bibinfo{title}{Relaxing bijectivity constraints with continuously
  indexed normalising flows},
\newblock in: \bibinfo{booktitle}{International conference on machine learning
  (2020) 2133--2143}, \bibinfo{organization}{PMLR}.
%Type = Article
\bibitem[{Zhang et~al.(2017)Zhang, Cisse, Dauphin, and
  Lopez-Paz}]{zhang2017mixup}
\bibinfo{author}{H.~Zhang}, \bibinfo{author}{M.~Cisse}, \bibinfo{author}{Y.~N.
  Dauphin}, \bibinfo{author}{D.~Lopez-Paz},
\newblock \bibinfo{title}{mixup: Beyond empirical risk minimization},
\newblock \bibinfo{journal}{arXiv preprint arXiv:1710.09412}
  (\bibinfo{year}{2017}).
%Type = Article
\bibitem[{Yang et~al.(2021)Yang, Meng, and Karniadakis}]{yang2021b}
\bibinfo{author}{L.~Yang}, \bibinfo{author}{X.~Meng}, \bibinfo{author}{G.~E.
  Karniadakis},
\newblock \bibinfo{title}{B-pinns: Bayesian physics-informed neural networks
  for forward and inverse pde problems with noisy data},
\newblock \bibinfo{journal}{Journal of Computational Physics}
  \bibinfo{volume}{425} (\bibinfo{year}{2021}) \bibinfo{pages}{109913}.

\end{thebibliography}

\newpage
%\appendix
\section*{\bf\LARGE{Supplementary Information}}

\renewcommand{\thefigure}{S\arabic{figure}}

In the Supplementary Information, for convenience of notation, we let $p(\cdot)$ denote the ``true'' density given by the data distribution and $\tilde p(\cdot)$ denote the density obtained from the augmented data generated by the GIN.
$q_{E}(z|x)$ and $q_{D}(y|z)$ denote the conditional densities defined by the encoder \eqref{eq:encoder} and the decoder.
Notice that for both training data and synthetic data, $p(z|x)$ and $\tilde p(z|x)$ are exactly known as
\[
p(z|x)=\tilde p(z|x)=q_E(z|x)
\]
for the given encoder. Hence,
\begin{eqnarray*}
p(x,y,z) & = & p(x,y)q_{E}(z|x),\\
\tilde{p}(\tilde x,\tilde z) & = & \tilde{p}(\tilde x)q_{E}(\tilde z|\tilde x).
\end{eqnarray*}
The decoder density $q_D(y|z)$ can be interpreted as a variational approximation of $p(y|z)$.

In experiments, both $q_E$ and $q_D$ are set to be Gaussian distributions, where $q_{E}(z|x)$ is a Gaussian distribution with mean $\mathrm{diag}(m(x))\bar z(x)$ and covariance $\mathrm{diag}(\mathbf{1}-m(x))^2$,  $q_{D}(y|z)$ is a Gaussian distribution with mean $\mu_{D}(z)$ and covariance $\Sigma_{D}(z)$.

\section*{\bf\large{\mycommand\label{app:vib} \quad Variational approximation of the information bottleneck objective}}

Here, we show the derivation of the variational IB objective \eqref{eq:vib_loss} for the sake of completeness, which is similar to the derivation in \cite{alemi2016deep}.

Let us start with the first term $I(Z;Y)$ in \eqref{eq:obj-aug}, i.e. the mutual information between the latent variable $Z$ and the output variables $Y$.
    Notice that
\begin{eqnarray*}
I(Z;Y) & = & \mathbb{E}_{p(y,z)}\left[\log\frac{p(y,z)}{p(y)p(z)}\right]\\
 & = & \mathbb{E}_{p(y,z)}\left[\log p(y|z)\right]-\mathbb{E}_{p(y)}\left[\log p(y)\right]\\
 & = & \mathbb{E}_{p(y,z)}\left[\log q_{D}(y|z)\right]+\mathbb{E}_{p(z)}\left[\mathbb{E}_{p(y|z)}\log\left[\frac{p(y|z)}{q_{D}(y|z)}\right]|z\right]-\mathbb{E}_{p(y)}\left[\log p(y)\right]\\
 & = & \mathbb{E}_{p(y,z)}\left[\log q_{D}(y|z)\right]+\mathbb{E}_{p(z)}\left[\mathrm{KL}\left(p(y|z)||q_{D}(y|z)\right)|z\right]+\mathcal{H}(Y),
\end{eqnarray*}
The second term on the r.h.s of the above equation is the mean value of the Kullback–Leibler (KL) divergence between $q_{D}(y|z)$ and $p(y|z)$, and the last term in the last inequality is the entropy of labels $\mathcal H(Y)$ that is independent of the optimization procedure and thus can be ignored during the training process. Due to the nonnegativity of the KL divergence, we have
\begin{equation}\label{eq:I_ZY_bound}
I(Z;Y)-\mathcal{H}(Y)\ge \widehat{I(Z;Y)}\triangleq\mathbb{E}_{p(y,z)}\left[\log q_{D}(y|z)\right],
\end{equation}
and the equality holds if $q_{D}(y|z)=p(y|z)$.

We next consider the term $I(\tilde Z;\tilde X)$, i.e. the mutual information between the augmented latent and input
    variables. Notice that
\begin{eqnarray}
I(\tilde Z;\tilde X) & = & \mathbb{E}_{\tilde{p}(\tilde x,\tilde z)}\left[\log q_{E}(\tilde z|\tilde x)\right]-\mathbb{E}_{\tilde{p}(\tilde z)}\left[\log\tilde{p}(\tilde z)\right]\nonumber\\
 & = & \mathbb{E}_{\tilde{p}(\tilde x,\tilde z)}\left[\log q_{E}(\tilde z|\tilde x)\right]-\mathbb{E}_{\tilde{p}(\tilde z)}\left[\log e(\tilde z)\right]-\mathrm{KL}\left(\tilde{p}(\tilde z)||e(\tilde z)\right)\nonumber\\
 & \le & \widehat{I(\tilde Z;\tilde X)}\label{eq:I_ZX_bound}
\end{eqnarray}
with
\begin{equation}\label{eq:I_ZX_bound_def}
\widehat{I(\tilde Z;\tilde X)}\triangleq\mathbb{E}_{\tilde{p}(\tilde x,\tilde z)}\left[\log q_{E}(\tilde z|\tilde x)\right]-\mathbb{E}_{\tilde{p}(\tilde z)}\left[\log e(\tilde z)\right]
\end{equation}
where $e$ is an arbitrary probability density function and the equality holds if $e(\tilde z)=\tilde p(\tilde z)$.

By combining \eqref{eq:I_ZY_bound} and \eqref{eq:I_ZX_bound}, we get
\[
\mathcal L_{\mathrm {VIB}} = \widehat{I(Z;Y)} - \beta\widehat{I(\tilde Z;\tilde X)}\le \mathcal L_{\mathrm {IB}}
-\mathcal H(Y)\]
and the equality can be attained with $q_D(y|z)=p(y|z), e(\tilde z)=\tilde p(\tilde z)$.
\iffalse
here $p(z|x)$ is designed by users and can be exactly known as a parametric model
$$
p(z \mid x)=q_{E}(z \mid x).
$$
In this paper, we model $q_{E}(z|x)$ as
$$
z=\operatorname{diag}(m(x)) \bar{z}(x)+\operatorname{diag}(1-m(x)) z_{0},
$$
where $m(x)$, $\bar{z}(x)$ are deep neural networks, each element of $m(x)$ is in [0,1], and $z_{0} \sim \mathcal{N}(0,I)$.
Here $q_{E}(z \mid x)$ is a trainable stochastic encoder transforms $X$ to $Z$, and $m(x)$ can be interpreted as the probability of $x$ being from the training data distribution.
$m(x)=1$ implies that $x$ is close to one of training data and the corresponding label $y$ can then be reliably predicted with the deterministic latent variable $z$.
If $m(x) = 0$, $x$ is considered as an out-of-distribution data and $z$ becomes an uninformative variable $z_{0}$.
In addition, $e(z)$ is a model of the variational marginal distribution of $z$ which is employed as Normalizing flow (see \ref{app:Realnvpintroduction}) in our IB-UQ model. 

Combining equation (\ref{IB_yz}) and (\ref{IB_zx}) provides a simple tractable lower bound for the IB objective,
\begin{equation}
\label{IB_lb}
\mathcal{L}_{IB}(Z)\geq \widehat{\mathcal{L}}_{IB}(Z):=\widehat{I(Y ; Z)}-\beta \widehat{I(Z ; X)},
\end{equation}
with 
\begin{equation}
\label{IB_lb1}
\widehat{I(Y ; Z)}:= \mathbb{E}_{p(y, z)}\left[\log q_{D}(y \mid z)\right],\quad \widehat{I(Z ; X)}:= \mathbb{E}_{p(x, z)}\left[\log q_{E}(z \mid x)\right]-\mathbb{E}_{p(z)}[\log e(z)].
\end{equation}
\fi

\section*{\bf\large{\mycommand\label{app:Realnvpintroduction} \quad Normalizing Flow}}
Let $W\in \R^D$ be a known random variable with tractable probability density function $p_W(w)$ (e.g. standard Gaussian distribution) and $V$ be an unknown random variable with density function $p_V(v)$. The normalizing flow model seeks to find a bijective mapping $f_I$ satisfying $V=f_I(W)$. Then we can compute the probability density function of the random variable $V$ by using the change of variables formula:
\begin{equation}\label{eqn:change_variable}
p_V(v)=\bigg |\text{det} \frac{\partial f_I^{-1}(v)}{\partial v} \bigg | p_W(w),
\end{equation}
where $f_I^{-1}$ is the inverse map of $f_I$. Assume the map $f_I$ is parameterized by $\theta$. Given a set of observed training data $\mathcal{D}=\{v_i\}_{i=1}^{N}$ from $V$, the log likelihood is given as following by taking log operator on both sides of (\ref{eqn:change_variable}):
\begin{equation}\label{eqn:likelihood1}
    \begin{aligned}
    \text{log}~p_V(\mathcal{D}|\theta) &=\sum_{i=1}^{N}\text{log}~p_V(v_i|\theta)\\
    & = \sum_{i=1}^{N}\text{log}~p_W(w_i)-\text{log}~\bigg |\text{det} \frac{\partial f_I(w_i|\theta)}{\partial w_i} \bigg |,
    \end{aligned}
\end{equation}
with $w_i=f_I^{-1}(v_i)$. The parameters $\theta$ can be optimized during the training by maximizing the above log-likelihood.

In this work, we will adopt the Non-Volume Preserving (RealNVP) model proposed by Dinh et al.~in ~\cite{dinh2016density} to construct the invertible transformation $f_I$, and $p_W$ is selected to be a standard Gaussian distribution.
\iffalse
The main idea of RealNVP is to construct $f$ by first split $z$ two parts, i.e. $w=(w_1,w_2)$ with $w_1\in \mathbb{R}^d$ and $w_2\in \mathbb{R}^{D-d}$ for $1<d<D$. Then do the following transformation
\begin{equation}\label{eqn:realnvpdinh}
\begin{aligned}
x_1& = z_1\\
x_2&= z_2\odot \text{exp}(s(z_1))+t(z_1)
\end{aligned}
\end{equation}
where $s$ and $t$ are scaling and translation depending on $z_1$ modeled by neural networks, $\odot$ is Hadamard product. By doing this, we can easily obtain the Jacobian of $f$ as follows
\begin{equation*}
\left[
  \begin{array}{cc}
    \mathbb{I}_d & \mathbf{0}_{d\times(D-d)} \\
    \frac{\partial x_2}{\partial z_1} & \text{diag}~(\text{exp}(s(z_1))) \\
  \end{array}
\right].
\end{equation*}
Given the observation that this Jacobian is triangular, we can efficiently compute its determinant used in the normalizing flow model. The inverse map $f^{-1}:x\to z$ can be computed directly 
\begin{equation}\label{eqn:realnvpdinh1}
\begin{aligned}
z_1& = x_1\\
z_2&= (x_2-t(x_1))\odot \text{exp}(-s(x_1))
\end{aligned}
\end{equation}
To guarantee both parts can be updated, we can change the position of $z_1$ and $z_2$ alternatively in the implementation and perform the transformation in a composite form
$(y_1,y_2)=f(z_1,z_2), (x_2,x_1)=f'(y_2,y_1)$, where $f'$ is a function of the same form as $f$.
The question remains of how to choose $d$. This is done by splitting the dimensions in half in our numerical examples.
\fi
For more details on this topic, readers are referred to \cite{kobyzev2020normalizing} and references therein.

\section*{\bf\large{\mycommand\label{app:gin} \quad General Incompressible-Flow Network and Mixup}}
The General Incompressible-Flow Networks (GIN) model is similar in form to RealNVP but remains the volume-preserving dynamics, i.e. $\bigg |\text{det} \frac{\partial f_I(w)}{\partial w} \bigg |\equiv 1$, which can be done by subtracting the mean of logarithm of outputs from each scaling layer as in \cite{sorrenson2019disentanglement}. 
For given samples of a random variable $V$, we can train the parameters of GIN by maximum likelihood to approximate the data distribution as
$p_V(v)=p_W(f_I^{-1}(v))$, and new samples of $V$ can be generated by $v=f(w)$ with $w\sim p_W(w)$,
where $p_W$ is a standard Gaussian distribution.
As shown in \cite{dibak2022temperature}, if we scale the variable $W$ by a scalar $\tau>0$ in GIN, we can obtain a new random variable $V'=f_I(\sqrt{\tau}W)$ with distribution $p_{V'}(v')\propto p_V(v')^{\frac{1}{\tau}}$.

In our experiments, we utilize a GIN to approximate the data distribution $p(x)$ and generate samples distributed according to $\tilde p(x)\propto p(x)^{\frac{1}{\tau}}$.
It is worth pointing out that normalizing flows usually over-estimate the density of areas between modes of multi-modal distributions \cite{cornish2020relaxing}. But this limitation is not critical for applications of the GIN in IB-UQ, since the over-estimation of density can lead to more OOD samples for data augmentation.

According to our numerical experience, if the size of the training data set is extremely small, the empirical data distribution based approximation of $\widehat{I(Z;Y)}$ may yield overfitting. We can optionally
smooth the empirical training data distribution by employing data augmentation technique such as Mixup \cite{zhang2017mixup},
and perform estimation based on the new $x$, $y$:

\begin{equation}
\label{mixup}
\begin{aligned}
    & x :=\lambda x+(1-\lambda) x^{\prime} \\
    &y :=\lambda y+(1-\lambda) y^{\prime}
\end{aligned}
\end{equation}
where the interpolation ratio $\lambda\in [0,1]$ is sampled from Beta($\alpha, \alpha$) with a small $\alpha$ which will be specified later in the numerical examples, and ($x$, $y$) and ($x^{\prime}, y^{\prime}$) are independently sampled from the training data. The interpolated samples, still denotes as $(x,y)$ for simplicity, are then used in the training process to compute $\widehat{I(Z ; Y)}$. (See Line 4 in Algorithm \ref{alg:nff}.)
\iffalse
Specifically, if we draw $(x_i,y_i)_{i=1}^{B}$ from the training data set and let $(x_b,y_b)_{b=1}^B$ is the mixup dataset after a random permutation of $(1,...,B)$, then $\widehat{I(Z ; Y)}$ can be estimated as 
\begin{equation}
\label{IB_lb1_appr}
\widehat{I(Z ; Y)}\approx \frac{1}{B}\sum_{b} \log q_{D}\left(y_{b} \mid z_{b}\right),
\end{equation}
here $z_b$ is latent random variable distributed according to $q_E(z_b|x_b)$.
\fi
In contrast to the original Mixup paper \cite{zhang2017mixup} where $\alpha$ is typically set to $0.1\sim 0.4$, we found that a much smaller value of $\alpha$ ($0.005$) can provide satisfactory uncertainty quantification results in our experiments.

\section*{\bf\large{\mycommand\label{app:algorithm} \quad Information bottleneck based UQ for DNN function regression algorithm}}
According to (\ref{eq:I_ZY_bound}) and (\ref{eq:I_ZX_bound_def}), we can use a mini-batch of samples $\{(x_1,y_1),\ldots,(x_B,y_B)\}$ from training dataset and $\{\tilde x_1,\ldots,\tilde x_B\}$ generated by GIN to obtain unbiased estimates of $\widehat{I(Z;Y)}$ and $\widehat{I(Z;X)}$ as
\begin{eqnarray*}
\widehat{I(Z;Y)} & \approx & \frac{1}{B}\sum_{b=1}^{B}\log q_{D}(y_{b}|z_{b}),\\
\widehat{I(\tilde{Z};\tilde{X})} & \approx & \frac{1}{B}\sum_{b=1}^{B}\log q_{E}(\tilde{z}_{b}|\tilde{x}_{b})-\log e(\tilde{z}_{b}),
\end{eqnarray*}
where $z_b, \tilde z_b$ are sampled from $q_E(\cdot|x_b), q_E(\cdot|\tilde x_b)$. Therefore, we can maximize $\mathcal L_{\mathrm{VIB}}$ by stochastic gradient denscent, and the proposed IB-UQ algorithm is summarized in Algorithm \ref{alg:nff}.

\begin{algorithm}[H]
\begin{itemize}
\item\textbf{1.} Given the training set $\{(x_i,y_i)\}_{i=1}^N$, hyperparameters $\beta, \tau, \alpha$, batch size $B$ and learning rate $\eta$
\item\textbf{2.} Utilize the GIN model to obtain $\tilde p(x)\propto {p(x)}^{\frac{1}{\tau}}$ based on $\{x_1,\ldots,x_N\}$
\item\textbf{3.}  Randomly draw a mini-batch $(x_{1}, y_{1}), ... , (x_{B}, y_{B})$ from the training data
\item\textbf{4.}  Perform Mixup (optional):\\\\
    \quad (a) Generate a random permutation $(I_{1}, ... , I_{B})$ of $(1, ... ,B)$, and sample $\lambda_{1},\ldots,\lambda_B$
    from Beta$(\alpha, \alpha)$ \\\\
    \quad (b) Let $x_{b}:=\lambda_{b} x_{b}+\left(1-\lambda_{b}\right) x_{I_{b}}$ and $y_{b}:=\lambda_{b} y_{b}+\left(1-\lambda_{b}\right) y_{I_{b}}$ for $b = 1,...,B$
\item\textbf{5.} Draw $\tilde x_1,\ldots, \tilde x_B \sim\tilde p(x)$ by GIN
\item\textbf{6.} Let
\begin{eqnarray*}
z_{b} & = & \mathrm{diag}\left(m\left(x_{b}\right)\right)\bar{z}\left(x_{b}\right)+\mathrm{diag}\left(1-m\left(x_{b}\right)\right)z_{0,b},\\
\tilde{z}_{b} & = & \mathrm{diag}\left(m\left(\tilde{x}_{b}\right)\right)\bar{z}\left(\tilde{x}_{b}\right)+\mathrm{diag}\left(1-m\left(\tilde{x}_{b}\right)\right)\tilde{z}_{0,b},
\end{eqnarray*}
where $z_{0, b}, \tilde z_{0, b} \sim \mathcal{N}(0,I)$ for $b = 1, ... ,B$
\item\textbf{7.} Estimate $\mathcal L_{\mathrm{VIB}}$ by
\[
\widehat{{\mathcal L}_{\mathrm{VIB}}}=\frac{1}{B}\sum_{b=1}^B \log q_D(y_b|z_b) - \beta \left(\log q_E(\tilde z_b|\tilde x_b)-\log e(\tilde z_b)\right)
\]
    \item\textbf{8.} Update all the parameters $\theta$ of $q_{D}, q_E, e$ as
    $$
    \theta:=\theta+\eta\frac{\partial}{\partial \theta}\widehat{{\mathcal L}_{\mathrm{VIB}}}
    $$
    \item\textbf{9.} Repeat Steps 3$\sim$7 until convergence
\end{itemize}
\caption{\label{alg:nff}Information bottleneck based UQ for function regression}
\end{algorithm}

In this paper, to model $\bar z(x)$ in the encoder, we employ an multilayer perceptron (MLP) with linear output activation, while $\log m(x)$ is modeled using an MLP with LogSigmoid output activation. To ensure that $m(x)$ approaches $0$ for low-density samples of $p_G(x)$, we modify the output as
\[
\log m(x)^{\mathrm{new}} = \mathrm{LogSigmoid}\left(\log\frac{m(x)^{-}}{1-m(x)^{-}}+\log\left(\mathrm{tanh}\left(\frac{\mathrm{relu}\left(\log \tilde p(x)-\log \tilde p^{\prime}\right)}{\log \tilde p^{\prime\prime}-\log \tilde p^{\prime}}\right)+10^{-12}\right)\right),
\]
where $m(x)^-=\min\{m(x),\mathrm{Sigmoid}(10^{6})\}$, $\tilde p^{\prime}, \tilde p^{\prime\prime}$ are the $1$st and $5$th percentiles of $p_G(x)$. It can be easily verified that $m(x)^{\mathrm{new}}\approx m(x)^-\approx m(x)$ if $\tilde p(x)\gg \tilde p^{\prime\prime}$ and $m(x)\lessapprox {10}^{-6}$ for $\tilde p(x)\le\tilde p^{\prime}$.
In the decoder, we assume that $\Sigma_D(z)$ is a diagonal matrix, and we use an MLP with linear output activation to model $\mu_D(z)$ and the logarithm of the diagonal elements of $\Sigma_D(z)$. Additionally, we model $e(z)$ using a normalizing flow. Detailed architectures of the neural network are described in \ref{app:h-param}.

After training, for a new input $x$, we can draw a large number of samples of $(Y,Z)$ from the distribution $q_E(z|x)\cdot q_D(y|z)$, and perform the prediction and UQ by calculating the mean and variance of samples of $Y$.

\section*{\bf\large{\mycommand\label{app:GIN-analysis} \quad Analysis of data augmentation with GIN}}

For simplicity of notation, we denote by $\mathcal M=(q_E,q_D,e)$ the set of all density models involved in IB-UQ, and denote by $\mathcal M^*=(q_E^*,q_D^*,e^*)$ the optimal model that maximizes the variational IB objective \eqref{eq:vib_loss}. Here we shall prove the following property of the optimal encoder: \textit{If there exists a constant $C>1$ so that $p(y)^{-1}p(y|x)\le C$ for all $x, y$, the Kullback-Leibler divergence between $q_D^*(\cdot|x)$ and $e(\cdot)$ is not larger than $\epsilon \beta^{-1}\log C$ for all $x$ and $\epsilon $ satisfying $p(x)\le \epsilon \tilde p(x)$.}

Here we define
\begin{eqnarray*}
L_{YZ}(x;\mathcal M) & = & \iint\frac{p(x)}{\tilde p(x)}p(y|x)q_{E}(z|x)\log \frac{q_D(y|z)}{p(y)}\mathrm{d}y\mathrm{d}z,\\
L_{XZ}(x;\mathcal M) & = & \int q_{E}(z|x)\log\frac{q_{E}(z|x)}{e(z)}\mathrm{d}z,\\
L(x;\mathcal M) & = & L_{YZ}(x;\mathcal M)-\beta L_{XZ}(x;\mathcal M)
\end{eqnarray*}
where $p(\cdot)$ denotes the density obtained from the data distribution, and $\tilde p(\cdot)$ denotes the density of inputs generated by the GIN. Then the variational IB objective can be written as
\[
\mathcal L_{\mathrm {VIB}}(\mathcal M) + \mathcal H(Y) = \int \tilde p(x)L(x;\mathcal M)\mathrm dx,
\]
which implies that the optimal $\mathcal M^*$ maiximizes the value of $L(x;\mathcal M)$ for all $x$. Furthermore, it can be seen that $L_{XZ}(x;\mathcal M)$ is equal to the the Kullback-Leibler divergence $\mathrm{KL}(q_E(\cdot|x)||e(\cdot))$ for a given $x$.

Considering $L(x;\mathcal M)=0$ if $q_D(y|z)=p(y)$ and $q_E(z|x)=e(z)$, we get
\begin{equation}\label{eq:L_lb}
L(x;\mathcal M^*)\ge 0.
\end{equation}
In addition, according to the assumption of $p(y)^{-1}p(y|x)\le C$ and the condition for equality in inequality \eqref{eq:I_ZY_bound}, we have
\begin{eqnarray*}
q_{D}^{*}(y|z) & = & \frac{\int p(y|x)q_{E}^{*}(z|x)p(x)\mathrm{d}x}{\int q_{E}^{*}(z|x')p(x')\mathrm{d}x'}\\
 & \le & \frac{\int Cp(y)q_{E}^{*}(z|x)p(x)\mathrm{d}x}{\int q_{E}^{*}(z|x')p(x')\mathrm{d}x'}\\
 & = & Cp(y)
\end{eqnarray*}
and
\begin{eqnarray}
L_{YZ}(x;\mathcal{M}^*) & \le & \log C\cdot\iint\frac{p(x)}{\tilde p(x)}p(y|x)q_{E}^*(z|x)\mathrm{d}y\mathrm{d}z\nonumber\\
 & \le & \epsilon\log C\label{eq:LYZ_ub}
\end{eqnarray}
if $p(x)\le \epsilon \tilde p(x)$. Combining inequalities (\ref{eq:L_lb}, \ref{eq:LYZ_ub}) leads to
\begin{eqnarray*}
\mathrm{KL}(q_{E}^{*}(\cdot|x)||e^{*}(\cdot)) & = & L_{XZ}(x;\mathcal{M}^{*})\\
 & \le & \beta^{-1}L_{YZ}(x;\mathcal{M}^{*})\\
 & \le & \epsilon\beta^{-1}\log C
\end{eqnarray*}

Thus, for an OOD $x$ satisfying $p(x)\le \epsilon \tilde p(x)$ with a small $\epsilon$, the latent variable $z$ is approximately independent of $x$ for the optimal encoder obtained by maximizing the data-augmentation based IB objective.

\section*{\bf\large{\mycommand\label{app:operator-algorithm} \quad Information bottleneck based UQ for operator learning}}

Similar to conclusion in \ref{app:vib}, we have
\begin{eqnarray*}
I((Y,Z);s(Y)) & = & \mathbb{E}_{p(s(y),y,z)}\left[\log\frac{p(s(y),y,z)}{p(s(y))p(y,z)}\right]\\
 & = & \mathbb{E}_{p(s(y),y,z)}\left[\log p(s(y)|y,z)\right]-\mathcal{H}(s(Y))\\
 & \ge & \mathbb{E}_{p(s(y),y,z)}\left[\log q_{D}(s(y)|y,z)\right]-\mathcal{H}(s(Y))
\end{eqnarray*}
and
\begin{eqnarray*}
I(\tilde Z;\tilde O) & = & \mathbb{E}_{\tilde{p}(\tilde o,\tilde z)}\left[\log q_{E}(\tilde z|\tilde o)\right]-\mathbb{E}_{\tilde{p}(\tilde z)}\left[\tilde{p}(\tilde z)\right]\\
 & \le & \mathbb{E}_{\tilde{p}(\tilde o,\tilde z)}\left[\log q_{E}(\tilde z|\tilde o)\right]-\mathbb{E}_{\tilde{p}(\tilde z)}\left[e(\tilde z)\right].
\end{eqnarray*}
Hence,
\[
\mathcal L_{\mathrm{VIBONet}}\le\mathcal L_{\mathrm{IBONet}}-\mathcal H(s(Y)),
\]
and the equality holds if $q_D(s(y)|y,z)=p(s(y)|y,z)$ and $e(\tilde z)=\tilde p(\tilde z)$.
We can perform stochastic gradient descent to maximize $\mathcal L_{\mathrm{VIBONet}}$ as shown in Algorithm \ref{alg:onet}.  

\section*{\bf\large{\mycommand\label{app:additional-results} \quad Additional results}}

{\bf{Discontinuous function regression} }

\setcounter{figure}{0}

\begin{figure}[!htp]
\centering
\includegraphics[width=0.42\linewidth]{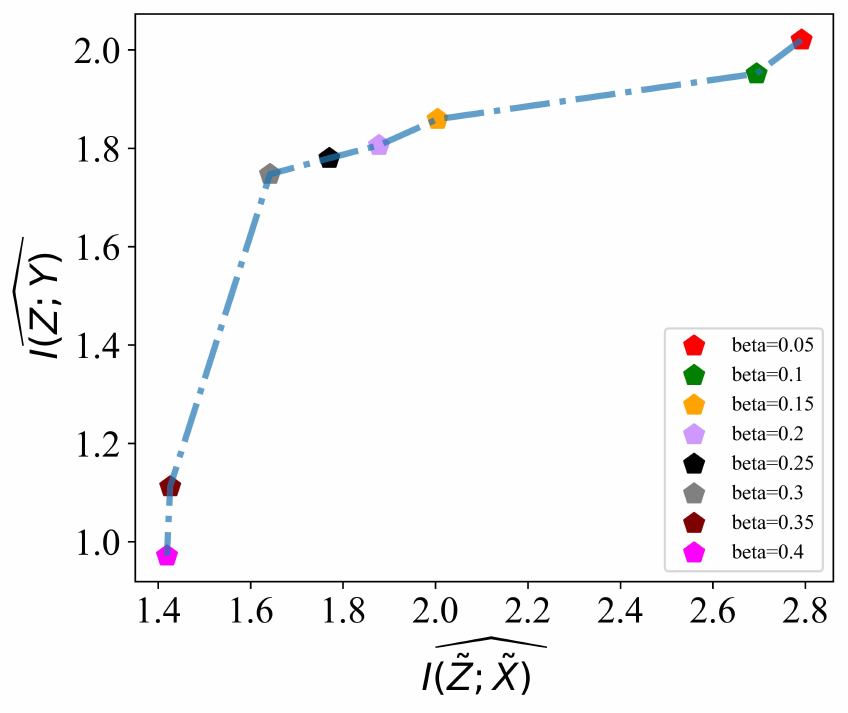}
\caption{Function regression problem (\ref{eq:DG function}) with noise scale $\sigma_u=0.1$. Visualization of the prediction and compression in the IB objective for different values of the IB Lagrange multiplier $\beta$ at the final training stage. } 
\label{fig:informationplane}
\end{figure}

We now conduct a study on how different values of the regularization parameter $\beta$ in the IB objective affect the model output. We use the information plane to visualize the performance of different $\beta$ in terms of the approximation of $\widehat{I(Z ; Y)}$ and  
$\widehat{I(\tilde Z ; \tilde X)}$ in the loss function at the final training stage. We run the IB-UQ code three times for every $\beta$, and pick up the one which has the smallest loss to report the corresponding values of $\widehat{I(Z ; Y)}$ and  
$\widehat{I(\tilde Z ; \tilde X)}$ in Figure \ref{fig:informationplane}. As the motivation of information bottleneck is to maximize the mutual information of $Y$ and $Z$, thus we can see from the plot that the optimal regularization parameter $\beta$ should vary in $[0.2, 0.3]$. Thus, for this function regression example, we choose $\beta=0.3$ to obtain the model prediction. The mixup parameter in (\ref{mixup}) is sampled from Beta($\alpha, \alpha$) with a small $\alpha=0.005$. 

\begin{algorithm}[H]
\begin{itemize}
\item\textbf{1.} Given the training set $\{(u^l,s^l)\}_{l=1}^N$, where measurements $o^l=\{u^l(x_j)\}_{j=1}^m$ and $\{s^l(y_j^l)\}_{j=1}^M$ are available, hyperparameters $\beta, \tau, \alpha$, batch sizes $B,B_M$ and learning rate $\eta$
\item\textbf{2.} Utilize the GIN model to obtain $\tilde p(o)\propto {p(o)}^{\frac{1}{\tau}}$ based on $\{o^1,\ldots,o^N\}$
\item\textbf{3.}  Randomly draw a mini-batch $(u^{1}, s^{1}), ... , (u^{B}, s^{B})$ from the training data, and sensor locations\footnotemark $\{y_1^l,\ldots,y_{B_M}^l\}$ for each $s^l$
\item\textbf{4.}  Perform Mixup (optional):\\\\
    \quad (a) Generate a random permutation $(I_{1}, ... , I_{B})$ of $(1, ... ,B)$, and sample $\lambda_1,\ldots,\lambda_B$
    from Beta$(\alpha, \alpha)$ \\\\
    \quad (b) Let $o^{b}:=\lambda_{b} o^{b}+\left(1-\lambda_{b}\right) o^{I_{b}}$ and $s^b(y_j^l):=\lambda_{b} s^b(y_j^l)+\left(1-\lambda_{b}\right) s^{I_b}(y_j^l)$ for $b = 1,...,B$ and $j=1,\ldots,B_M$
\item\textbf{5.} Draw $\tilde o^1,\ldots, \tilde o^B \sim\tilde p(o)$ by GIN
\item\textbf{6.} Let
\begin{eqnarray*}
z^{b} & = & \mathrm{diag}\left(m\left(o^{b}\right)\right)\bar{z}\left(o^{b}\right)+\mathrm{diag}\left(1-m\left(o^{b}\right)\right)z_0^b,\\
\tilde z^{b} & = & \mathrm{diag}\left(m\left(\tilde o^{b}\right)\right)\bar{z}\left(\tilde o^{b}\right)+\mathrm{diag}\left(1-m\left(\tilde o^{b}\right)\right)\tilde z_0^b,
\end{eqnarray*}
where $z_0^b, \tilde z_0^b \sim \mathcal{N}(0,I)$ for $b = 1, ... ,B$
\item\textbf{7.} Estimate $\mathcal L_{\mathrm{VIBONet}}$ by
\begin{eqnarray*}
\widehat{\mathcal{L}_{\mathrm{VIBONet}}} & = & \frac{1}{BB_M}\sum_{b=1}^{B}\sum_{j=1}^{B_M}\log q_{D}(s(y_{j}^{b})|y_{j}^{b},z^{b})\\
 &  & -\frac{\beta}{B}\sum_{b=1}^{B}\left(\log q_{E}(\tilde{z}^{b}|\tilde{o}^{b})-\log e(\tilde{o}^{b})\right)
\end{eqnarray*}
    \item\textbf{8.} Update all the parameters $\theta$ of $q_{D}, q_E, e$ as
    $$
    \theta:=\theta+\eta\frac{\partial}{\partial \theta}\widehat{{\mathcal L}_{\mathrm{VIBONet}}}
    $$
    \item\textbf{9.} Repeat Steps 3$\sim$8 until convergence
\end{itemize}
\caption{\label{alg:onet}Information bottleneck based UQ for operator learning}
\end{algorithm}
\footnotetext{In our experiments, sensor locations are fixed and we draw $y_1,\ldots,y_{B_M}$ for all $s^l$ in each iteration.}

\begin{figure}[htp!]
\centering
\centering
\includegraphics[width=0.75\linewidth]{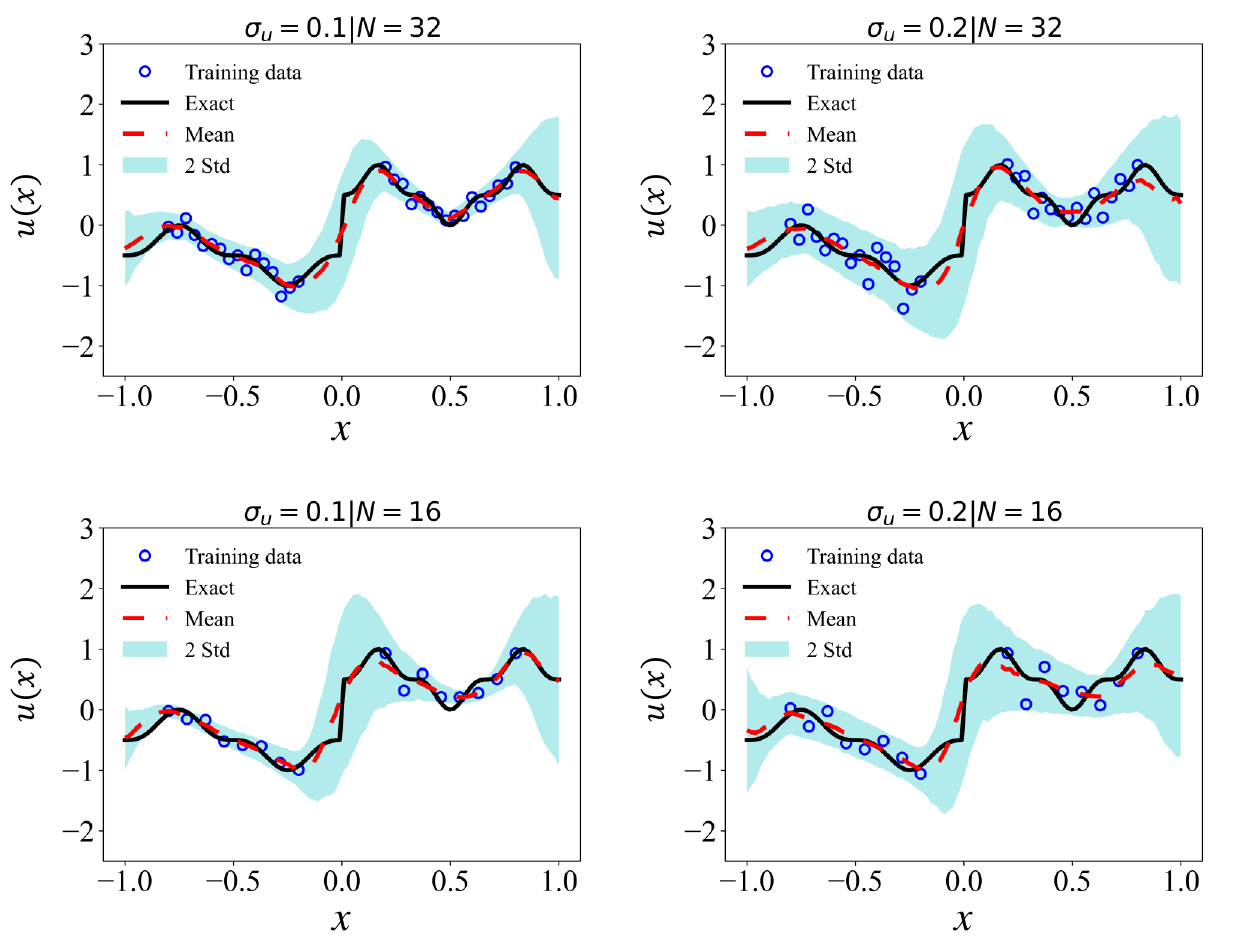}
\caption{Function regression problem (\ref{eq:DG function}). 
The uncertainty of IB-UQ increases with increasing noise magnitude and decreasing data set size. This is illustrated in the presented results, which depict the training data, exact function, and mean and uncertainty (standard deviations) obtained by the IB-UQ method for two noise scales (left: $\sigma_u = 0.1$; right: $\sigma_u = 0.2$) and two data set sizes (top: $N = 32$; bottom: $N = 16$).}
\label{fig:fun_2}
\end{figure}
We further test the performance of the IB-UQ method for cases with larger noise scales and smaller training data set sizes. Here, the parameter of mixup $\alpha$ is taken as 0.005 and 0.01 for training data set size 32 and 16 respectively. As shown in Figure \ref{fig:fun_2}, the standard
deviations (i.e., uncertainty) is expected to increase with increasing data noise scale and
decreasing data set size. 

{\bf{Climate model}}

\begin{figure}[htp!]
\centering
\centering
\subfloat[]{\includegraphics[width=0.99\linewidth]{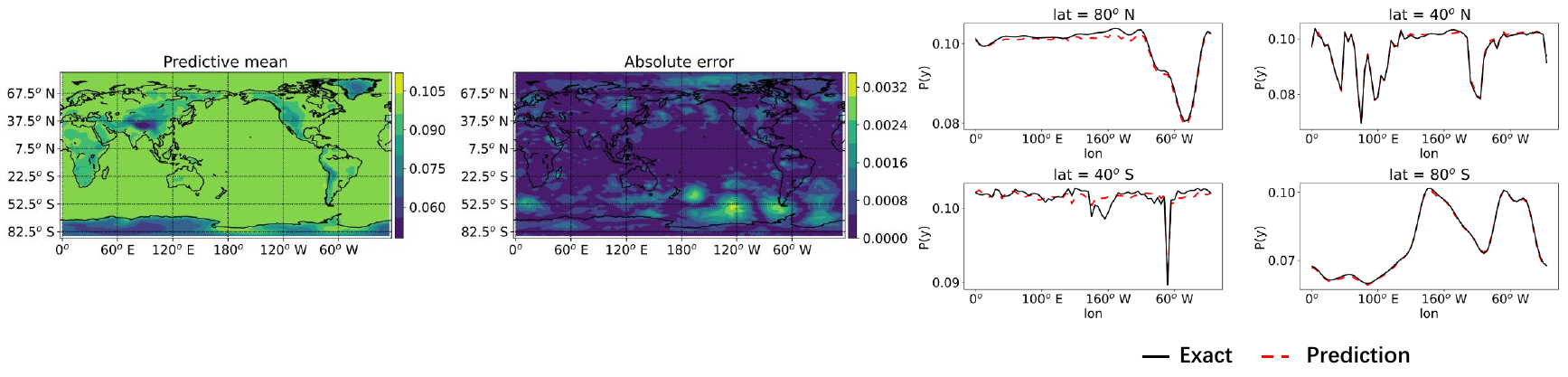}}\\
\subfloat[]{\includegraphics[width=0.99\linewidth]{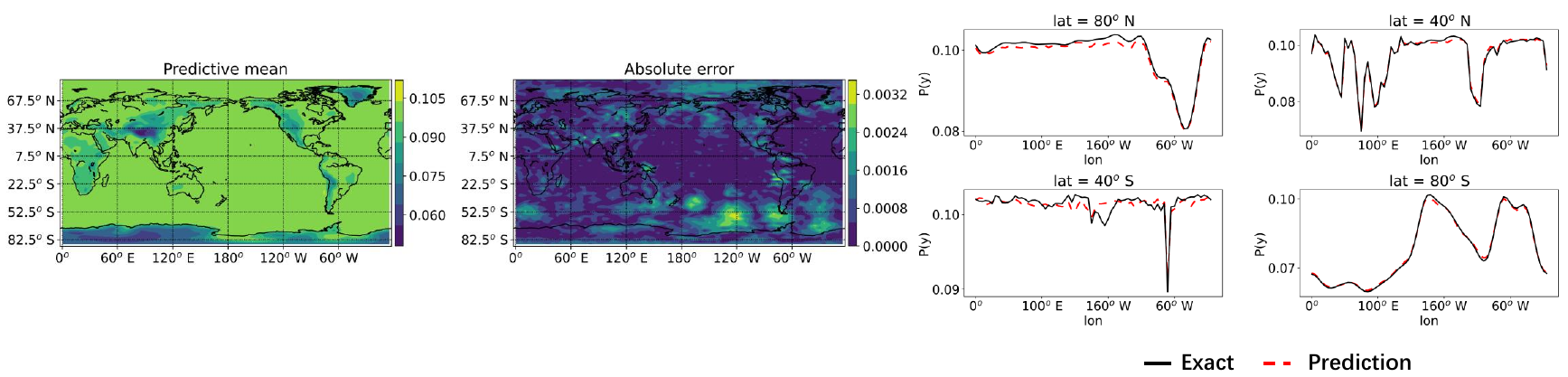}}
\\
\subfloat[]{\includegraphics[width=0.99\linewidth]{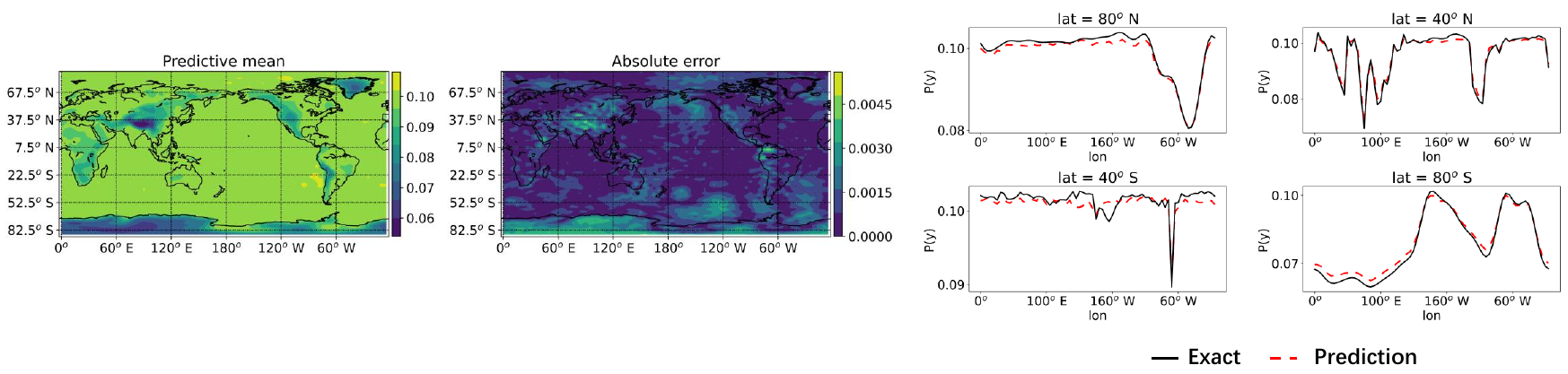}}\\
\caption{Representative results of pressure prediction by vanilla DeepONet algorithm for given surface air temperature. Data processing and other experiment details in Figure (a)-(c) are consistent with those in Figure \ref{fig:climate_ib}. Compared with the results given by our IB-UQ method in \ref{fig:climate_ib}, we can see vanilla DeepONet method can not provide the model's confidence without any uncertainty estimates.}
\label{fig:deeponet _climate_model}
\end{figure}

In the large-scale climate model, we used the vanilla DeepONet method to learn the mapping from surface temperature to surface pressure. Figure \ref{fig:deeponet _climate_model} illustrates the results of this experiment.  This case demonstrates the capability of our proposed IB-UQ model to provide meaningful uncertainty estimates simultaneously for operator learning, which is an improvement over the vanilla DeepONet method.

\section*{\bf\large{\mycommand\label{app:h-param} \quad Hyperparameters and neural network architectures used in the numerical examples}}

In this supplementary Information, we provide all hyperparameters and neural network architectures, as well as  hyperparameters used in data augmentation for IB-UQ method in all the numerical experiments. For the function approximation problem, Gaussian process regression that corresponds to BNN with $\mathcal{N}(0,1)$ as prior distribution is considered as in \cite{yang2021b}. The parameters for HMC method and Deep Ensemble method are inherited from \cite{psaros2022uncertainty} (see Table \ref{Hyperparameter-func}). The prior for $\theta$ (the weights and bias in the neural networks) is set as independent standard Gaussian distribution for each component. In HMC, the mass matrix is set to the identity matrix, i.e., $\mathbf{M} = \mathbf{I}$, the leapfrog step is set to $L = 50\epsilon$, the time step is $\epsilon$ = 0.1, the burn-in steps are set to 2000, and the total number of samples is 1000. For the Deep Ensemble methods, we run the DNN approximation code for $M=20$ times to get the ensemble results. The details of neural networks architecture settings of RealNVP, GIN, $m(x)$, $\bar z(x)$ and the decoder net used in IB-UQ methods for function regression are summarized in Table \ref{Hyperparameter-func and operator}. The temperatures $\tau$ of the GIN model are set to 16 (discontinuous function regression) and 4 (California housing prices prediction) and the latent dimension is 20. The Adam optimizer with batch size $B=256$ is used to optimize all models with a base learning rate $10^{-3}$, which we decrease by a factor of 0.1. The batch size of sensor locations in operator learning is $B_M=101$ (see Line 3 in Algorithm \ref{alg:onet}), i.e., for each output function in the mini-batch, we select measurements at $101$ query locations points randomly drawn from the grid points for training. For the discontinuous function regression problem, considering the training data size is small, we repeated the available data for each mini-batch in order to maintain a batch size of $256$, which resulted in the generation of more noisy data due to Mixup.

For learning the nonlinear operator of partial differential equations using the proposed IB-UQ model, we parametrize the branch and trunk networks in the DeepONet method using a multilayer perceptron with 3 hidden layers and 128 neurons per layer. The details of neural networks architecture settings of RealNVP, GIN, $m(x)$, $\bar z(x)$ and the decoder net used in IB-UQ model are summarized in Table \ref{Hyperparameter-func and operator}. The temperature $\tau$ in the GIN model is set to 1.4 and the latent dimension is 16. The Adam optimizer with batch size 256 is used to optimize all models with a base learning rate $10^{-3}$, which we decrease by a factor of 0.1.

Additionally, as part of the training procedure, we adopt a pre-training phase and set $\beta=0$ during the initial epochs.
According to our experience, this strategy significantly enhances the stability and overall performance of the training results.
All the implementation settings for the training algorithm are summarized in Table \ref{Hyperparameter-func_and_onet}.

\begin{table}[h!]\small
  \begin{center}
\caption{Hyperparameter settings: HMC and Deep Ensemble method used in the function approximation problems.}
    \label{Hyperparameter-func}
   \begin{tabular}{c|c|c|c|c|c}
     \hline
      Method & \multicolumn{5}{c}{Hyperparameter}\\
     \hline
      \multirow{2}*{HMC} & $\epsilon$ & $T$ & burn-in steps & $M$  & net-size\\
      \cline{2-6}
       & $0.1$ & $50$ & $2 \times 10^{3}$ & $10^{3}$ & $50 \times 2$\\
      \hline
      \multirow{2}*{Deep Ensemble} & train-steps & $M$ & weight-decay & init & net-size \\
      \cline{2-6}
       & $2 \times 10^{3}$ & 20 & $4 \times 10^{-3}$ & Xavier (normal) & $50 \times 2$\\
      \hline
    \end{tabular} 
  \end{center}
\end{table}

\begin{table}[h!]\footnotesize
  \begin{center}
    \caption{Hyper-parameter setting: IB-UQ architecture used in the function approximation  and operator learning problems.}
    \label{Hyperparameter-func and operator}
    \setlength{\tabcolsep}{0.6mm}
    \begin{tabular}{c|c|c|c|c|c|c}
    \hline
    & $m(x)$ & $\bar{z}(x)$ & RealNVP & GIN-RealNVP &  \multicolumn{2}{c}{\multirow{2}*{}} \\
    & width$\times$ depth  & width$\times$ depth & blocks | width$\times$ depth & blocks | width$\times$ depth &  \multicolumn{2}{c}{\multirow{2}*{}}\\
    \hline
    \multirow{3}*{Function regression} &  &  &   &   &   \multicolumn{2}{c}{Decoder-Nets} \\
    &  \normalsize{32 $\times$ 2} & \normalsize{32 $\times$ 2}  & \normalsize{3 \quad 256 $\times$ 6}  & \normalsize{3 \quad 256 $\times$ 6}  &  \multicolumn{2}{c}{width$\times$ depth} \\
    \cline{6-7}
    &  &  &  &  &  \multicolumn{2}{c}{32 $\times$ 2} \\
    \hline
    \multirow{3}*{Operator learning} &  &  &  &  &  branch-net & trunk-net \\
    & \normalsize{128 $\times$ 3} & \normalsize{128 $\times$ 3} & \normalsize{3 \quad 256 $\times$ 6} & \normalsize{3 \quad 256 $\times$ 6} &  width$\times$ depth & width$\times$ depth \\
    \cline{6-7}
    &  &  &  &  &  128 $\times$ 3 &  128 $\times$ 3\\
    \hline
    \end{tabular} 
  \end{center}
\end{table}

\begin{table}[h!]\small
  \begin{center}
    \caption{Hyper-parameter settings: Model training settings used in the function approximation and operator learning problems.}
    \label{Hyperparameter-func_and_onet}
    \setlength{\tabcolsep}{2mm}
    \begin{tabular}{c|c|c|c|c|c}
     \hline
     & Hyperparameter & Learning rate & Optimizer & Epochs & Learning rate decay\\
     \hline
     \multirow{2}*{Function regression} & GIN &  \multirow{4}*{\large{${10}^{-3}$}} & \multirow{4}*{\large{Adam}} & 200 & 0.1/50 epochs  \\
     \cline{2-2} \cline{5-6}
     & IB-UQ &  &     & \makecell[c]{5000 (including 1000 \\ pre-training epochs)} & 0.1/2000 epochs \\
     \cline{1-2} \cline{5-6}
     \multirow{2}*{Operator learning} & GIN &  &     & 100 &  0.1/50 epochs  \\
     \cline{2-2} \cline{5-6}
     & IBUQONet &  &  & \makecell[c]{600 (including 100 \\ pre-training epochs)} & 0.1/200 epochs \\
     \hline
    \end{tabular} 
  \end{center}
\end{table}

\end{document}